
\input amstex
\UseAMSsymbols
\documentstyle {amsppt}
\loadbold \loadeufm

\def\trace{\operatorname{trace}}
\def\sp{\operatorname{Span}}
\def\ricm{\widehat{\operatorname{Ric}}}
\def\Aut{\operatorname{Aut}}
\def\pt{\frac{\partial}{\partial t}}
\def\pu{\frac{\partial}{\partial u}}
\def\pv{\frac{\partial}{\partial v}}
\def\a{\alpha}
\def\b{\beta}
\def\g{\gamma}
\def\t{\tau}
\def\l{\lambda}
\def\tl{\tilde{\lambda}}
\def\th{\theta}
\def\f{\varphi}
\def\R{\Bbb R}

\def\v{\bold v}

\pagewidth{37.5 pc} \pageheight{55 pc} 

\document
\topmatter
 \title Affine hypersurfaces admitting a pointwise symmetry
\endtitle 
\author Y. Lu, C. Scharlach \endauthor
\thanks The second author was partially supported by the DFG-project 
{\it Geometric Problems and Special PDEs}.
\endthanks

\address C. Scharlach, Technical University Berlin, Fak. II, 
Institut of Mathematics, D-10623 Berlin, Germany
\endaddress
\address Y. Lu, School of Mathematics, Xiamen University, 361000, 
Xiamen, Fujian, China
\endaddress



\abstract An affine hypersurface $M$ is said to admit a
pointwise symmetry, if there exists a subgroup $G$ of $\Aut(T_p M)$
for all $p\in M$, which preserves (pointwise) the affine metric $h$,
the difference tensor $K$ and the affine shape operator $S$. In this
paper, we deal with positive definite affine hypersurfaces of
dimension three. First we solve an algebraic problem. We determine the
non-trivial stabilizers $G$ of the pair $(K,S)$ under the action of
$SO(3)$ on an Euclidean vectorspace $(V,h)$ and find a representative
(canonical form of $K$ and $S$) of each $SO(3)/G$-orbit. Then, we
classify hypersurfaces admitting a pointwise $G$-symmetry for all
non-trivial stabilizers $G$ (apart of $Z_2$). Besides well-known
hypersurfaces (for $Z_2\times Z_2$ we get the locally homogenous
hypersurface $(x_1-1/2\, x_3^2)(x_2-1/2\, x_4^2)=1$) we obtain
e.g. warped product structures of two-dimensional affine spheres
(resp. quadrics) and curves. \endabstract

\endtopmatter

\head 1. Introduction\endhead

 In this paper we study nondegenerate (equi-)affine hypersurfaces
$F\colon M^n \to \R^{n+1}$. It is well known that there exists a
canonical choice of transversal vector field $\xi$ called the affine
(Blaschke) normal, which induces a connection $\nabla$, a symmetric
bilinear form $h$ and a 1-1 tensor field $S$ by
$$\align
&D_X Y =\nabla_X Y +h(X,Y)\xi,\tag 1.1\\
&D_X \xi =-SX,\tag 1.2
\endalign$$
for all $X,Y \in {\Cal X}(M)$. The connection $\nabla$ is called
the induced affine connection, $h$ is called the affine metric (or
Blaschke metric) and $S$ is called the affine shape operator.  In
general $\nabla$ is not the Levi Civita connection
$\widehat\nabla$ of $h$. The difference tensor $K$ is defined as
$$K(X,Y)=\nabla_X Y-\widehat\nabla_X Y,\tag 1.3$$ for all $X,Y \in
{\Cal X}(M)$. Moreover the form $h(K(X,Y),Z)$ is a symmetric cubic
form with the property that for any fixed $X\in {\Cal X}(M)$, $\trace
K_X$ vanishes.  This last property is called the apolarity condition.
The difference tensor $K$, together with the affine metric $h$ and the
affine shape operator are the most fundamental algebraic invariants
for a nondegenerate affine hypersurface. We say that $M$ is positive
definite if the affine metric $h$ is positive definite. For details of
the basic theory of nondegenerate affine hypersurfaces we refer to
\cite{6} and
\cite{9}.

 Furthermore, the fundamental equations for an affine hypersurface
 are,
 $$R(X,Y)Z=h(Y,Z)SX-h(X,Z)SY,\tag1.4$$
 $$\nabla h(X,Y,Z)=\nabla h(Y,X,Z),\tag1.5$$
 $$\nabla S(X,Y)=\nabla S(Y,X),\tag1.6$$
 $$h(SX,Y)=h(X,SY),\tag1.7$$
 where $R$ is the curvature tensor with respect to $\nabla$.
 Define $C:=\nabla h$, it is totally symmetric and
 $$C(X,Y,Z)=-2h(K(X,Y),Z).\tag1.8$$

An affine hypersurface $M$ is said to admit a pointwise $G$-symmetry,
if $G$ is a (orientation-preserving) subgroup of $\Aut(T_p M)$ for all
$p\in M$, which preserves (pointwise) $h$, $K$ and $S$. Necessarily
$G$ is a subgroup of the isometry group of $T_p M$.

The study of submanifolds which admit a pointwise symmetry was
initiated by Bryant in \cite{1} where he studied $3$-dimensional
special Lagrangian submanifolds of $\Bbb C^3$. Here the symmetry group
preserves (pointwise) the first fundamental form and the fundamental
cubic. Following essentially the same approach, a classification of
$3$-dimensional positive definite affine hyperspheres admitting
pointwise symmetries was obtained in
\cite{11}. Here we deal with all $3$-dimensional positive
definite hypersurfaces admitting a pointwise symmetry, i.~e. the
symmetry group preserves additionally to the affine metric and the
difference tensor the affine shape operator.

  First we solve an algebraic problem. In Section 2 we determine the
non-trivial stabilizers $G$ of the pair $(K,S)$ under the action of
$SO(3)$ on a Euclidean vectorspace $(V,h)$. The non-trivial
stabilizers are isomorphic to a copy of $SO(3)$, $Z_2\times SO(2)$,
$SO(2)$, $A_4$, $S_3$, $Z_2\times Z_2$, $Z_3$ or $Z_2$. Compared with
the stabilizers of $K$ (cf. \cite{1}) we get additionally $Z_2\times
SO(2)$ and $Z_2\times Z_2$. Furthermore we find a representative
(canonical form of $K$ and $S$) of each $SO(3)/G$-orbit.

In Section 3 we start with the classificaton of positive definite
affine hypersurfaces admitting a pointwise $G$-symmetry for all
non-trivial stabilizers $G$ (apart of $Z_2$). There are no
hypersurfaces admitting a pointwise $Z_2\times SO(2)$-symmetry.  We
will show that a hypersurface with $G=S_3$ must be a hypersphere,
hence the classification can be found in \cite{11}. A hypersurface
admits a pointwise $Z_2\times Z_2$-symmetry if and only if it is
affine equivalent to $(x_1-1/2\, x_3^2)(x_2-1/2\, x_4^2)=1$ (and thus it
is a locally homogeneous affine hypersurface with rank one shape
operator).

The final section is organized as follows. First, we show that for
$G=Z_3$ resp. $G=SO(2)$ we can extend the canonical form of $K$ and
$S$ locally and thus obtain information about the coefficients of $K$,
$S$ and $\nabla$ from the basic equations of Gauss, Codazzi and
Ricci. In particular, it follows that the hypersurface admits a warped
product structure $\R\times_{e^f}N^2$. Then following essentially the
same approach as in
\cite {12}, we classify such hypersurfaces by showing how they can be
constructed starting from $2$-dimensional positive definite affine
spheres resp. quadrics. 

This last classification can be seen as a generalization of the well
known Calabi product of hyperbolic affine spheres and of the
constructions for affine spheres considered in
\cite{4}. The following natural question for a (de)composition theorem, 
related to the Calabi product and its generalizations in
\cite{4}, gives another motivation for studying $3$-dimensional
hypersurfaces admitting a pointwise symmetry:

\proclaim {(De)composition Problem }Let $M^n$ be a nondegenerate affine
hypersurface in $\R^{n+1}$. Under what conditions do there exist
  affine hyperpsheres $M_1^r$ in $\R^{r+1}$ and $M_2^s$ in $\R^{s+1}$,
  with $r+s=n-1$, such that $M = I \times_{f_1} M_1 \times_{f_2} M_2$,
  where $I \subset \R$ and $f_1$ and $f_2$ depend only on $I$
  (i.e. $M$ admits a warped product structure)? How can the original
  immersion be recovered starting from the immersion of the affine
  spheres?
\endproclaim

Of course the first dimension in which the above problem can be
considered is three and our study of $3$-dimensional affine
hypersurfaces with $Z_3$-symmetry or $SO(2)$-symmetry provides an
answer in that case.

\head 2. Canonical Forms for $K$ and $S$  \endhead 

Let $(V,h)$ be a 3-dimensional Euclidean space, endowed with a
symmetric endomorphism $S$ and a symmetric $(2,1)$-form $K$ with
$\trace K_X=0$ for all $X\in V$. Let $G$ be a subgroup of the isometry
group $SO(3)$ of $(V,h)$. Assume that $S$ and $K$ are preserved by
$G$, in other words, for any $X,Y\in V, g\in G\subset SO(3)$, we have
 $$h(gX,gY)=h(X,Y),\tag I$$
 $$K(gX,gY)=g(K(X,Y)),\tag II$$
 $$S(gX)=g(SX).\tag III$$

 If $\{e_1,e_2,e_3\}$ is an orthonormal basis,
 we may define a polynomial invariant to $G$, that is
 $$f(x,y,z)=h(K(xe_1+ye_2+ze_3,xe_1+ye_2+ze_3),xe_1+ye_2+ze_3).$$
 According to \cite{1} (resp. \cite{11}), we know that
 \proclaim{ Theorem 2.1} There exists an
 orthonormal basis of $V$ such that either
 \smallskip
 $(i)$ $f=0$, in this case $f$ is preserved by $SO(3)$,
\smallskip
 $(ii)$ $f=\l(2x^3-3xy^2-3xz^2)$, for some $\l>0$ in which case $f$ is
 preserved by $SO(2)$, which consists of any rotation on $yz-$plane,
\smallskip
 $(iii)$ $f=6\l xyz$, for some $\l>0$ in which case $f$ is preserved
 by $A_4$ of order 12, which is generated by rotation by an angle of
 $\pi$ about $x-,y-,z-$axes, and rotation by an angle of
 $\frac{2\pi}{3}$ about the line $x=y=z$, 
\smallskip 
$(iv)$ $f=\l(x^3-3xy^2)$ for some $\l>0$, in which case $f$ is
 preserved by $S_3$ of order 6, which is generated by rotation by an
 angle of $\pi$ about $x-$axis, and rotation by an angle of
 $\frac{2\pi}{3}$ about the $z-$axis, 
\smallskip 
$(v)$ $f=\l(2x^3-3xy^2-3xz^2)+6\mu xyz$, for some $\l,\mu>0$, with
 $\l\neq\mu$, in which case $f$ is preserved by $Z_2$ of order 2,
 which is generated by rotation by an angle of $\pi$ about $x-$axis,
 \smallskip $(vi)$ $f=\l(2x^3-3xy^2-3xz^2)+\mu(y^3-3yz^2)$ for some
 $\l,\mu>0$, with $\mu\neq\sqrt{2}\l$, in which case $f$ is preserved
 by $Z_3$, which is generated by rotation by an angle of
 $\frac{2\pi}{3}$ about the $x-$axis.\endproclaim

 \smallskip
 
 It is obvious that the traceless $(2,1)$-form $K$ can be determined
 completely by $f$. As to the polynomials in Theorem 2.1, they are
 respectively
 \smallskip
 (i) $K=0$;
\smallskip
 (ii) $K_{e_1}=\pmatrix 2\l&0&0\\0&-\l&0\\0&0&-\l\endpmatrix$,
 $K_{e_2}=\pmatrix 0&-\l&0\\-\l&0&0\\0&0&0\endpmatrix$,
 $K_{e_3}=\pmatrix 0&0&-\l\\0&0&0\\-\l&0&0\endpmatrix$;
\smallskip
 (iii) $K_{e_1}=\pmatrix 0&0&0\\0&0&\l\\0&\l&0\endpmatrix$,
 $K_{e_2}=\pmatrix 0&0&\l\\0&0&0\\\l&0&0\endpmatrix$,
 $K_{e_3}=\pmatrix 0&\l&0\\\l&0&0\\0&0&0\endpmatrix$;
 \smallskip
 (iv) $K_{e_1}=\pmatrix \l&0&0\\0&-\l&0\\0&0&0\endpmatrix$,
 $K_{e_2}=\pmatrix 0&-\l&0\\-\l&0&0\\0&0&0\endpmatrix$,
 $K_{e_3}=0$;
 \smallskip
 (v) $K_{e_1}=\pmatrix 2\l&0&0\\0&-\l&\mu\\0&\mu&-\l\endpmatrix$,
 $K_{e_2}=\pmatrix 0&-\l&\mu\\-\l&0&0\\\mu&0&0\endpmatrix$,
 $K_{e_3}=\pmatrix 0&\mu&-\l\\\mu&0&0\\-\l&0&0\endpmatrix$;
 \smallskip
 (vi) $K_{e_1}=\pmatrix 2\l&0&0\\0&-\l&0\\0&0&-\l\endpmatrix$,
 $K_{e_2}=\pmatrix 0&-\l&0\\-\l&\mu&0\\0&0&-\mu\endpmatrix$,
 $K_{e_3}=\pmatrix 0&0&-\l\\0&0&-\mu\\-\l&-\mu&0\endpmatrix$;
 \smallskip
 \smallskip
 Note that an orthogonal transformation $T$ can create a new basis
 $\{\tilde{e}\}=\{\tilde{e}_1,\tilde{e}_2,\tilde{e}_3\}$ from the
 given one, $\{{e}\}=\{{e}_1,{e}_2,{e}_3\}$. The relation of
 these two bases is determined by,
 $$\{\tilde{e}_1,\tilde{e}_2,\tilde{e}_3\}=\{{e}_1,{e}_2,{e}_3\}T_e,$$
 where $T_e$ is the matrix form of $T$ with respect to the basis
 $\{e\}$. In this case, the matrix form of any linear transformation $L$ with
 respect to $\{\tilde{e}\}$ comes from a similar
 transformation $T$ acting on the matrix with respect to $\{e\}$,
 i.e.
 $$L_{\tilde{e}}=T_e^{-1}L_eT_e.$$
 
Therefore, in $\{\tilde{e}\}$ we obtain a congruent expression for the
 group $G$.  In general, the expressions of $f$ are different in these
 two bases. However, they will be the same when $T$ is chosen from the
 groups mentioned in Theorem 2.1, and then condition (I) and (II) are
 satisfied. We would like to find out for which subgroups also $S$ is
 preserved, i.~e.  in matrix form, $SG=GS$. (By Theorem 2.1 the group
 $G$ is a subgroup of $SO(3)$, $SO(2)$, $A_4$, $S_3$, $Z_2$ or $Z_3$.)
 In the following, we divide the possibilities of $G$ into three
 cases. First, $G$ is finite and cyclic; second, it is finite but
 non-cyclic; finally, it is infinite.
\smallskip
 From now on we shall denote by $P_i$ the rotation by an angle of
 $\pi$ about $e_i$, $R_i$ the rotation by an angle of
 $\frac{2\pi}{3}$ about $e_i$, $Q$ the rotation by an angle of
 $\frac{2\pi}{3}$ about $x=y=z$.
 \smallskip
 For the sake of proof, we present here the structure of
 $S_3$ and $A_4$,
 $$\align S_3&=\langle P_1,R_3\rangle\\
 &=\{id\}\cup\{P_1,P_1R_3,R_3P_1\}\cup\{R_3,R_3^2\},\tag 2.1\\
 A_4&=\langle P_1,P_2,Q\rangle\\
 &=\{id\}\cup\{P_1,P_2,P_3\}\cup\{Q,Q^2,P_1Q,P_2Q,P_3Q,Q^2P_1,Q^2P_2,Q^2P_3\},\tag
 2.2
 \endalign$$
 in which the second part of the union is the set of elements
 whose order is $2$, and the third part is of order $3$. Note that 
$$\align P_1 S= S P_1 \quad
&\Rightarrow \quad S=\pmatrix a&0&0\\0&b&d\\0&d&c\endpmatrix,\\
R_3 S= S R_3 \quad
&\Rightarrow \quad S=\pmatrix a&0&0\\0&a&0\\0&0&b\endpmatrix,\\
Q S= S Q \quad
&\Rightarrow \quad S=\pmatrix a&b&b\\b&a&b\\b&b&a\endpmatrix.
\endalign$$

 {\bf Case 1.} Let $G$ be cyclic, i.e., $G=\langle
 g\rangle$, in which $g$ has finite order. Then we obtain the
 following lemmas.

 \proclaim{Lemma 2.1} Let $G=\langle g\rangle$ be a cyclic subgroup of
 $SO(3)$ and let $g$ be of order $2$. We have some orthonormal basis,
 such that
 $$K_{e_1}=\pmatrix 2\l&0&0\\0&-\l&\mu\\0&\mu&-\l\endpmatrix,
 K_{e_2}=\pmatrix 0&-\l&\mu\\-\l&0&0\\\mu&0&0\endpmatrix,
 K_{e_3}=\pmatrix 0&\mu&-\l\\\mu&0&0\\-\l&0&0\endpmatrix,$$
 $$S=\pmatrix a&0&0\\0&b&d\\0&d&c\endpmatrix,$$ where either $\l>0$,
 $\mu>0$, $\l\neq\mu$, or $\l=\mu>0$ and $(b-c)^2+(b+c-2a-2d)^2\neq0$,
 or $\l=0$, $\mu>0$ and $d\neq0$, or $\mu=0$, $\l>0$ and
 $(b-c)^2+d^2\neq0$. $K$ and $S$ are preserved by $G=Z_2=\langle
 P_1\rangle$.  \endproclaim
\demo{Proof}
 \roster
 \item $G=Z_2$: From Theorem 2.1 we know $G=\langle P_1\rangle$.
 Then the expression of $S$ comes from $P_1S=SP_1$, and $K$ from
 the theorem. In this case, $\l>0$, $\mu>0$, and $\l\neq\mu$;
 \item $G\subset S_3$: If $G=\langle P_1R_3\rangle$, then it is
 congruent to $\langle P_1\rangle$ since
 $(R_3P_1)^{-1}P_1R_3(R_3P_1)=P_1$. This indicates we can
 construct a new basis by $R_3P_1$, such that $G=\langle P_1\rangle$
 and $f$ is invariant. The same is true for $G=\langle R_3P_1\rangle$
 and $P_1R_3$. Thus we can assume that $G=\langle
 P_1\rangle$. Through the orthogonal transformation
 $$(\tilde{e}_1,\tilde{e}_2,\tilde{e}_3)=(e_1,e_2,e_3) \pmatrix
 1&&\\&\frac{1}{\sqrt{2}}&-\frac{1}{\sqrt{2}}
 \\&\frac{1}{\sqrt{2}}&\frac{1}{\sqrt{2}}\endpmatrix,\tag 2.3$$ $f=
 \tl(3x^3 -3 x y^2)$ can be turned into the form
 $\frac{\tl}{2}(2\tilde{x}^3-3\tilde{x}\tilde{y}^2-3\tilde{x}\tilde{z}^2)+\frac{\tl}{2}
 6\tilde{x}\tilde{y}\tilde{z}$, and for $S$ we get in matrix form
$$
 \pmatrix 1&&\\&\frac{1}{\sqrt{2}}&-\frac{1}{\sqrt{2}}
 \\&\frac{1}{\sqrt{2}}&\frac{1}{\sqrt{2}}\endpmatrix \pmatrix
 a&&\\&b&d\\&d&c\endpmatrix \pmatrix
 1&&\\&\frac{1}{\sqrt{2}}&\frac{1}{\sqrt{2}}
 \\&-\frac{1}{\sqrt{2}}&\frac{1}{\sqrt{2}}\endpmatrix =\pmatrix
 a&&\\&\frac12 (b-2d+c)&\frac12 (b-c)\\&\frac12 (b-c)&\frac12
 (b+2d+c)\endpmatrix,
$$
 Since $R_3\not\in G$, we have that $R_3 S\neq S R_3$, which implies that
 $b\neq c$ or $b+c\neq 2(a+d)$.  So we get the second possibility of the Lemma
 ($\l=\mu=\frac{\tl}{2}>0$); 
\item $G\subset
 A_4$: Since $Q^2P_2Q=QP_3Q^2=P_1$, we can assume $G=\langle
 P_1\rangle$. $f=6\tl xyz$ has the right form with
 $\l=0$ and $\mu=\tl>0$. To obtain a cyclic group, $S$ must not be
 invariant to any other element of $A_4$, therefore we have $d\neq0$;
 \item $G\subset SO(2)$: Clearly this is the case when $\mu=0$,
 $\l=\tl>0$. If there is a non-trivial element of $SO(2)$ preserving $S$,
 which is not $P_1$, then 
$$\pmatrix \cos\th&-\sin\th\\ \sin\th&\cos\th \endpmatrix 
\pmatrix b&d\\d&c\endpmatrix = \pmatrix b&d\\d&c\endpmatrix 
\pmatrix \cos\th&-\sin\th\\ \sin\th&\cos\th \endpmatrix \quad \text{with} \quad \th\neq k\pi$$ 
i.~e. $d=0$ and $b=c$. Since $G$ is cyclic, we have
 $(b-c)^2+d^2\neq0$; 
\item $G\subset SO(3)$: Due to $f=0$, we can
 choose an orthonormal basis such that $S$ is diagonal, in which case
 $S$ is invariant to any $P_i$. This implies that $G$ is not cyclic.
 \qed\endroster \enddemo

 \proclaim{Lemma 2.2} Let $G=\langle g\rangle$ be a cyclic subgroup of
 $SO(3)$ and let $g$ be of order $3$. We have some orthonormal basis,
 such that
 $$K_{e_1}=\pmatrix 2\l&0&0\\0&-\l&0\\0&0&-\l\endpmatrix,
 K_{e_2}=\pmatrix 0&-\l&0\\-\l&\mu&0\\0&0&-\mu\endpmatrix,
 K_{e_3}=\pmatrix 0&0&-\l\\0&0&-\mu\\-\l&-\mu&0\endpmatrix,$$
 $$S=\pmatrix a&0&0\\0&b&0\\0&0&b\endpmatrix,$$ where either $\l>0$,
 $\mu>0$ and $\mu\neq\sqrt{2}\l$, or $\mu=\sqrt{2}\l>0$ and $a\neq
 b$. $K$ and $S$ are preserved by $G=Z_3=\langle R_1\rangle$.
 \endproclaim
\demo{Proof}
 \roster
 \item $G=Z_3$: It comes simply from the theorem, and $\l>0$,
 $\mu>0$, $\mu\neq\sqrt{2}\l$;
 \item $G\subset S_3$: $G$ must be $\langle R_3\rangle$. So in
 this basis,
 $$
 S=\pmatrix a&&\\&a&\\&&b\endpmatrix.
 $$
 It is easy to check that $S$ is invariant to $P_1$, which
 indicates it is actually $S_3$-symmetry;
 \item $G\subset A_4$: If $G=\langle Q\rangle$, then the
 transformation
 $$(\tilde{e}_1,\tilde{e}_2,\tilde{e}_3)=(e_1,e_2,e_3) \pmatrix
 \frac{1}{\sqrt{3}}&\frac{2}{\sqrt{6}}&0\\
 \frac{1}{\sqrt{3}}&-\frac{1}{\sqrt{6}}&\frac{1}{\sqrt{2}}\\
 \frac{1}{\sqrt{3}}&-\frac{1}{\sqrt{6}}&-\frac{1}{\sqrt{2}}\endpmatrix$$
 yields $G=\langle R_1\rangle$ and $f=6\tl x y z$ is transformed into
 $\frac{\tl}{\sqrt{3}}(2\tilde{x}^3-3\tilde{x}\tilde{y}^2
 -3\tilde{x}\tilde{z}^2 +\sqrt{2} \tilde{y}^3- 3\sqrt{2}
 \tilde{y}\tilde{z}^2)$. We get the result when
 $\mu=\sqrt{2}\l=\sqrt{2}\frac{\tl}{\sqrt{3}}$. Since
 $(P_iQ)^2=Q^2P_i$, it suffices to show that $\langle P_iQ\rangle$ is
 congruent to $\langle Q\rangle$ in $A_4$. Taking $P_1Q$ as an
 example, one can check $(Q^2P_2)(P_1Q)(P_2Q)=Q$. The condition
 $R_1S=SR_1$ yields the expression of $S$ immediately. Assume $a=b$,
 then the symmetry is completely determined by $K$.  In this case $G$
 can be extended to $A_4$. Hence, in order to obtain a cyclic group,
 $a$ must not be equal to $b$; 
\item $G\subset SO(2)$: Obviously this is
 the case when $\l>0$ and $\mu=0$. Note that at this moment, $V$
 actually admits $SO(2)$-symmetry;
\item $G\subset SO(3)$: It is easy to see this is the result when
 $\l=\mu=0$, and $V$ also admits $SO(2)$-symmetry.
 \qed\endroster
 \enddemo

 If $G=\langle g\rangle$ is a cyclic subgroup of
 $SO(3)$ and $g$ is of order other than $2$ and $3$,
 then actually $SO(2)\subseteq G$, which is not finite.
 
\bigskip
 {\bf Case 2.} Let $G$ be finite but have several generators. Thus $G$
 only can be a subgroup of $S_3$, $A_4$, $SO(2)$ or $SO(3)$.  We get
 the following canonical forms with isotropy group $G$:
 
\proclaim{Lemma 2.3} Let $G$ be a finite, non-cyclic subgroup of
 $SO(3)$. There exists an orthonormal basis of $V$ such that either

\smallskip
(i) 
 $$K_{e_1}=\pmatrix 0&0&0\\0&0&\l\\0&\l&0\endpmatrix,
 K_{e_2}=\pmatrix 0&0&\l\\0&0&0\\\l&0&0\endpmatrix,
 K_{e_3}=\pmatrix 0&\l&0\\\l&0&0\\0&0&0\endpmatrix,$$
 $$S=\pmatrix a&0&0\\0&b&0\\0&0&c\endpmatrix,$$ in which $a$, $b$, $c$ are
 not all equal when $\l>0$, and distinct when $\l=0$. $K$ and $S$ are
 preserved by $G=Z_2\times Z_2=\langle P_1,P_2\rangle$; or 

\smallskip
 (ii) 
 $$K_{e_1}=\pmatrix \l&0&0\\0&-\l&0\\0&0&0\endpmatrix,
 K_{e_2}=\pmatrix 0&-\l&0\\-\l&0&0\\0&0&0\endpmatrix,
 K_{e_3}=0,$$
 $$S=\pmatrix a&0&0\\0&a&0\\0&0&b\endpmatrix,$$ where $ \l>0$. $K$ and
$S$ are preserved by $G=S_3$; or 

\smallskip (iii) 
 $$K_{e_1}=\pmatrix 0&0&0\\0&0&\l\\0&\l&0\endpmatrix,
 K_{e_2}=\pmatrix 0&0&\l\\0&0&0\\\l&0&0\endpmatrix,
 K_{e_3}=\pmatrix 0&\l&0\\\l&0&0\\0&0&0\endpmatrix,$$
 $$S=aI,$$ where $\l>0$.  $K$ and $S$ are preserved by $G=A_4$
 \endproclaim
 \demo{Proof} If $G\subseteq S_3$, then from (2.1) it must be $\langle P_1,
 R_3\rangle$, i.e., $G=S_3$. This corresponds to case (ii) of the lemma.

 Suppose $G\subseteq A_4$. If $G$ has two generators of order
 $3$, then $V$ has two transversal planes invariant to $S$ and thus
 $S=aI$, in which case $G=A_4$. If $G$ has one generator of order $2$ and another of
 order $3$, by the following relations,
 $$
 \langle P_i,P_jQ\rangle=\langle P_iP_jQ,P_jQ\rangle
 =\cases
 \langle Q,P_iQ\rangle\quad i=j\\
 \langle P_kQ,P_jQ\rangle \quad \{i,j,k\}=\{1,2,3\}
 \endcases
 $$
 we know $S=aI$ and $G=A_4$ again. Obviously this corresponds to case
 (iii) of the lemma.  If $G$ is generated by two elements of order
 $2$, we get $G=\langle P_i,P_j\rangle_{i\neq
 j}=\{id,P_1,P_2,P_3\}$. This corresponds to case (i) of the lemma
 with $\l>0$. Of course in this case $a$, $b$, $c$ are not all equal,
 otherwise $G$ could be extended to $A_4$.

 If $G\subset SO(2)$, then at least one generator is not of order
 $2$. So $S$ has an invariant plane, which must be the plane of
 $SO(2)$-rotation. Therefore, $G$ can be
 extended to $SO(2)$.

 If $G\subset SO(3)$, we have $K=0$ and the symmetry is
 completely determined by $S$. With respect to some orthonormal
 basis $S$ is diagonal. It is easy to see that $G$ is finite if
 and only if $S$ has three distinct eigenvalues, which corresponds to 
 case (i) of the lemma with $\l=0$.
 \qed
 \enddemo

\bigskip
{\bf Case 3.}  Let $G$ be infinite or continuous. Then we have
 $SO(2)\subseteq G\subseteq SO(3)$ and we get the following canonical
 forms with isotropy group $G$: 
\proclaim{Lemma 2.4} 
Let $G$ be an infinite subgroup of $SO(3)$. There exists an
 orthonormal basis of $V$ such that either \smallskip (i)
 $$K_{e_1}=\pmatrix 2\l&&\\&-\l&\\&&-\l\endpmatrix, K_{e_2}=\pmatrix
 0&-\l&0\\-\l&0&0\\0&0&0\endpmatrix, K_{e_3}=\pmatrix
 0&0&-\l\\0&0&0\\-\l&0&0\endpmatrix,$$
 $$S=\pmatrix a&0&0\\0&b&0\\0&0&b\endpmatrix,$$ 
where $\l>0$. $K$ and $S$ are preserved by $G=SO(2)$; or
 \smallskip
 (ii)  
 $$K=0,\quad S=\pmatrix a&0&0\\0&b&0\\0&0&b\endpmatrix,$$
where $a\neq b$. $K$ and $S$ are preserved by $G=Z_2\times SO(2)$; or
 \smallskip
 (iii)  
 $$K=0,\quad S=aI.$$ $K$ and $S$ are preserved by $G=SO(3)$.
 \endproclaim
 \demo{Proof} If $G=SO(2)$, then we have case (i) immediately.
 Otherwise $G\subseteq SO(3)$ and $K=0$. With respect to some
 orthonormal basis $S$ is diagonal. This gives the cases (ii)
 and (iii) immediately.
 \qed
 \enddemo

 As a conclusion we have 
\proclaim {Theorem 2.2} 
Let $(V,h)$ be a 3-dimensional Euclidean space, endowed with a
 symmetric endomorphism $S$ and a symmetric $(2,1)$-form $K$ with
 $\trace K_X=0$ for all $X\in V$. Assume that $K$ and $S$ are
 preserved by a subgroup $G$ of the isometry group $SO(3)$ of
 $(V,h)$. Then $G$ is isomorphic to either $SO(3)$, $Z_2\times SO(2)$,
 $SO(2)$, $A_4$, $S_3$, $Z_2\times Z_2$, $Z_3$ or $Z_2$.  With respect
 to some suitable orthonormal basis in $V$, we get the following
 canonical forms for $K$ and $S$:
\smallskip
(i) $$K=0,\quad S=aI,$$ in this case $K$ and $S$ are preserved by $SO(3)$;
\smallskip
(ii) $$K=0,\quad S=\pmatrix a&0&0\\0&b&0\\0&0&b\endpmatrix,\quad a\neq
b,$$ in this case $K$ and $S$ are preserved by $Z_2\times SO(2)$,
which consists of any rotation on $yz$-plane and $P_2$;
\smallskip
(iii)
 $$K_{e_1}=\pmatrix 2\l&0&0\\0&-\l&0\\0&0&-\l\endpmatrix,
 K_{e_2}=\pmatrix 0&-\l&0\\-\l&0&0\\0&0&0\endpmatrix,
 K_{e_3}=\pmatrix 0&0&-\l\\0&0&0\\-\l&0&0\endpmatrix,$$
 $$S=\pmatrix a&0&0\\0&b&0\\0&0&b\endpmatrix,\quad \l>0,$$ in this
case $K$ and $S$ are preserved by $SO(2)$, which consists of any
rotation on $yz$-plane, 
\smallskip
(iv)
 $$K_{e_1}=\pmatrix 0&0&0\\0&0&\l\\0&\l&0\endpmatrix,
 K_{e_2}=\pmatrix 0&0&\l\\0&0&0\\\l&0&0\endpmatrix,
 K_{e_3}=\pmatrix 0&\l&0\\\l&0&0\\0&0&0\endpmatrix,$$
 $$S=aI,\quad\l>0,$$ 
in this case $K$ and $S$ are preserved by $A_4=\langle P_1, P_2,Q\rangle$;
\smallskip
(v)
 $$K_{e_1}=\pmatrix \l&0&0\\0&-\l&0\\0&0&0\endpmatrix,
 K_{e_2}=\pmatrix 0&-\l&0\\-\l&0&0\\0&0&0\endpmatrix,
 K_{e_3}=0,$$
 $$S=\pmatrix a&0&0\\0&a&0\\0&0&b\endpmatrix,\quad\l>0,$$
in this case $K$ and $S$ are preserved by $S_3=\langle P_1,R_3\rangle$;
\smallskip
(vi)
 $$K_{e_1}=\pmatrix 0&0&0\\0&0&\l\\0&\l&0\endpmatrix,
 K_{e_2}=\pmatrix 0&0&\l\\0&0&0\\\l&0&0\endpmatrix,
 K_{e_3}=\pmatrix 0&\l&0\\\l&0&0\\0&0&0\endpmatrix,$$
 $$S=\pmatrix a&0&0\\0&b&0\\0&0&c\endpmatrix,$$ in which either $a$, $b$, $c$
 are not all equal when $\l>0$ or $a$, $b$, $c$ are distinct when
 $\l=0$; in this case $K$ and $S$ are preserved by $Z_2\times
 Z_2=\langle P_1, P_2\rangle$;
\smallskip
(vii)
 $$K_{e_1}=\pmatrix 2\l&0&0\\0&-\l&0\\0&0&-\l\endpmatrix,
 K_{e_2}=\pmatrix 0&-\l&0\\-\l&\mu&0\\0&0&-\mu\endpmatrix,
 K_{e_3}=\pmatrix 0&0&-\l\\0&0&-\mu\\-\l&-\mu&0\endpmatrix,$$
 $$S=\pmatrix a&0&0\\0&b&0\\0&0&b\endpmatrix,$$
 in which either $\l>0$, $\mu>0$ and $\mu\neq\sqrt{2}\l$, or
 $\mu=\sqrt{2}\l>0$ and $a\neq b$;
in this case $K$ and $S$ are preserved by $Z_3=\langle R_1\rangle$;
\smallskip
(viii)
 $$K_{e_1}=\pmatrix 2\l&0&0\\0&-\l&\mu\\0&\mu&-\l\endpmatrix,
 K_{e_2}=\pmatrix 0&-\l&\mu\\-\l&0&0\\\mu&0&0\endpmatrix,
 K_{e_3}=\pmatrix 0&\mu&-\l\\\mu&0&0\\-\l&0&0\endpmatrix,$$
 $$S=\pmatrix a&0&0\\0&b&d\\0&d&c\endpmatrix,$$
 in which either $\l>0$, $\mu>0$ and $\l\neq\mu$, or $\l=\mu>0$ and
 $(b-c)^2+(b+c-2a-2d)^2\neq0$, or $\l=0$, $\mu>0$ and $d\neq0$, or
 $\mu=0$, $\l>0$ and $(b-c)^2+d^2\neq0$;
in this case $K$ and $S$ are preserved by $Z_2=\langle P_1\rangle$.
\endproclaim

According to Theorem 2.2 we have 
\proclaim {Corollary 2.3} 
Let $p\in M^3$ and assume that there exists an orientation preserving
 isometry which preserves $K$ and $S$ at $p$. Then there exists an
 orthonormal basis of $T_pM^3$ such that either:
\smallskip
(i) $$K=0,\quad S=aI,$$ in this case $K$ and $S$ are preserved by $SO(3)$;
\smallskip
(ii) $$K=0,\quad S=\pmatrix a&0&0\\0&b&0\\0&0&b\endpmatrix,\quad a\neq
b,$$ in this case $K$ and $S$ are preserved by $Z_2\times SO(2)$,
which consists of any rotation on $yz$-plane and $P_2$;
\smallskip
(iii)
 $$K_{e_1}=\pmatrix 2\l&0&0\\0&-\l&0\\0&0&-\l\endpmatrix,
 K_{e_2}=\pmatrix 0&-\l&0\\-\l&0&0\\0&0&0\endpmatrix,
 K_{e_3}=\pmatrix 0&0&-\l\\0&0&0\\-\l&0&0\endpmatrix,$$
 $$S=\pmatrix a&0&0\\0&b&0\\0&0&b\endpmatrix,\quad \l>0,$$ in this
case $K$ and $S$ are preserved by $SO(2)$, which consists of any
rotation on $yz$-plane, 
\smallskip
(iv)
 $$K_{e_1}=\pmatrix 0&0&0\\0&0&\l\\0&\l&0\endpmatrix,
 K_{e_2}=\pmatrix 0&0&\l\\0&0&0\\\l&0&0\endpmatrix,
 K_{e_3}=\pmatrix 0&\l&0\\\l&0&0\\0&0&0\endpmatrix,$$
 $$S=aI,\quad\l>0,$$ 
in this case $K$ and $S$ are preserved by $A_4=\langle P_1, P_2,Q\rangle$;
\smallskip
(v)
 $$K_{e_1}=\pmatrix \l&0&0\\0&-\l&0\\0&0&0\endpmatrix,
 K_{e_2}=\pmatrix 0&-\l&0\\-\l&0&0\\0&0&0\endpmatrix,
 K_{e_3}=0,$$
 $$S=\pmatrix a&0&0\\0&a&0\\0&0&b\endpmatrix,\quad\l>0,$$
in this case $K$ and $S$ are preserved by $S_3=\langle P_1,R_3\rangle$;
\smallskip
(vi)
 $$K_{e_1}=\pmatrix 0&0&0\\0&0&\l\\0&\l&0\endpmatrix,
 K_{e_2}=\pmatrix 0&0&\l\\0&0&0\\\l&0&0\endpmatrix,
 K_{e_3}=\pmatrix 0&\l&0\\\l&0&0\\0&0&0\endpmatrix,$$
 $$S=\pmatrix a&0&0\\0&b&0\\0&0&c\endpmatrix,$$ in which either $a$,
 $b$, $c$ are not all equal when $\l>0$ or $a$, $b$, $c$ are distinct
 when $\l=0$; in this case $K$ and $S$ are preserved by $Z_2\times
 Z_2=\langle P_1, P_2\rangle$;
\smallskip
(vii)
 $$K_{e_1}=\pmatrix 2\l&0&0\\0&-\l&0\\0&0&-\l\endpmatrix,
 K_{e_2}=\pmatrix 0&-\l&0\\-\l&\mu&0\\0&0&-\mu\endpmatrix,
 K_{e_3}=\pmatrix 0&0&-\l\\0&0&-\mu\\-\l&-\mu&0\endpmatrix,$$
 $$S=\pmatrix a&0&0\\0&b&0\\0&0&b\endpmatrix,$$
 in which either $\l>0$, $\mu>0$ and $\mu\neq\sqrt{2}\l$, or
 $\mu=\sqrt{2}\l>0$ and $a\neq b$;
in this case $K$ and $S$ are preserved by $Z_3=\langle R_1\rangle$;
\smallskip
(viii)
 $$K_{e_1}=\pmatrix 2\l&0&0\\0&-\l&\mu\\0&\mu&-\l\endpmatrix,
 K_{e_2}=\pmatrix 0&-\l&\mu\\-\l&0&0\\\mu&0&0\endpmatrix,
 K_{e_3}=\pmatrix 0&\mu&-\l\\\mu&0&0\\-\l&0&0\endpmatrix,$$
 $$S=\pmatrix a&0&0\\0&b&d\\0&d&c\endpmatrix,$$
 in which either $\l>0$, $\mu>0$ and $\l\neq\mu$, or $\l=\mu>0$ and
 $(b-c)^2+(b+c-2a-2d)^2\neq0$, or $\l=0$, $\mu>0$ and $d\neq0$, or
 $\mu=0$, $\l>0$ and $(b-c)^2+d^2\neq0$;
in this case $K$ and $S$ are preserved by $Z_2=\langle P_1\rangle$.
\endproclaim

 \head 3. Hypersurfaces admitting $S_3$- or $Z_2\times Z_2$-symmetry \endhead

 In the following, we will assume that $M^3$ admits a pointwise
 $G$-symmetry, i.~e. that at every point $h$, $K$ and $S$
 are preserved by the group $G$. 

By Corollary 2.3 we can not have a pointwise $Z_2\times
 SO(2)$-symmetry, since the vanishing of $K$ implies that $M^3$ is a
 quadric and thus $S=a I$. In \cite{11}, a complete classification of
 three-dimensional positive definite affine hyperspheres admitting a
 pointwise $G$-symmetry was obtained. Since pointwise $Z_2$-symmetry
 is a rather weak assumption, in this case the author just got a
 completely integrable system, which in general can not be solved
 explicitly. So we will omit this case in the following
 classification.

Hence, due to Corollary 2.3, in the further we will only study the
 cases of pointwise $SO(2)$-, $S_3$-, $Z_2\times Z_2$- and
 $Z_3$-symmetry. In the rest of this section we will show that
 a hypersurface $M^3$ with $G=S_3$ must be a hypersphere, and then give a
 classification for $G=Z_2\times Z_2$. 

 A short computation shows that, for $G=S_3$, $K$ and $S$ are
 commutable, i.~e. $K(SX,Y)=K(X,SY)=SK(X,Y)$. So $\widehat{\nabla}S$
 is symmetric and consequently $M^3$ is a hypersphere ($h$ is positive
 definite) (cf. \cite{7}).

 Now let $G=Z_2\times Z_2=\langle P_1,P_2\rangle$, therefore at every
 point, we can find an orthonormal basis $\{e_1,e_2,e_3\}$ such that
 $K$ and $S$ are in canonical form (Case (vi) in Corollary 2.3). We
 want to extend the basis locally.

\proclaim {Lemma 3.1} 
If $M^3$ admits a pointwise $Z_2\times Z_2$-symmetry, then there exists
locally an orthonormal frame $\{e_1,e_2,e_3\}$ such that
 $$K_{e_1}=\pmatrix 0&0&0\\0&0&\l\\0&\l&0\endpmatrix, K_{e_2}=\pmatrix
 0&0&\l\\0&0&0\\\l&0&0\endpmatrix, K_{e_3}=\pmatrix
 0&\l&0\\\l&0&0\\0&0&0\endpmatrix,$$
 $$S=\pmatrix a&&\\&b&\\&&c\endpmatrix,$$ in which either $a$, $b$,
$c$ are not all equal when $\l>0$ or $a$, $b$, $c$ are distinct when
$\l=0$.
\endproclaim
\demo{Proof}
Let $p\in M^3$. If $S$ has three distinct eigenvalues, then there
 exists an open neighborhood $U_p$, in which we can select $\{e_i\}$
 as the smooth unit eigenvector fields, and $S$ has at any point of $U_p$
 three distinct eigenvalues. Moreover, difference
 tensor and shape operator are of the form above.

 If, in a neighborhood of $p$, $S$ has only two distinct eigenvalues,
 we first define the 1-dimensional unit normal eigenvector field
 as $e_3$. Let $F_3=e_3$, and $F_1,F_2$ be any smooth orthonormal vector
 fields in the 2-dimensional eigendistribution. At any point
 $q$ in such neighborhood, we have $\{e_i\}|_q$ defined pointwisely as
 above, and 
 $$\align &F_1=\cos\t e_1+\sin\t e_2,\\
 &F_2=-\sin\t e_1+\cos\t e_2,\endalign$$
 for some $\t(q)$. Thus, with respect to $\{F_i\}$
 $$K_{F_1}=\pmatrix
 0&0&\sin2\t\l\\0&0&\cos2\t\l\\\sin2\t\l&\cos2\t\l&0\endpmatrix$$
 Since $\l> 0$, we can define a smooth function $\t$ via
 $$\sin 2\t\cdot h(K_{F_1}F_2,F_3)=\cos 2\t\cdot h(K_{F_1}F_1,F_3),$$
and the orthonormal frame  
 $$\align &e_1=\cos\t F_1-\sin\t F_2,\\
 &e_2=\sin\t F_1+\cos\t F_2,\\&e_3=F_3.\endalign $$
 will have the properties we were looking for.  \qed
\enddemo

 Let $\widehat{\nabla}_{e_i}e_j=\f_{ij}^ke_k$ be the Levi-Civita
 connection. Obviously, $\f_{ij}^k=-\f_{ik}^j$, and $\f_{ij}^j=0$. It
 is well known, that we obtain the following Codazzi equations for $C$
 from the fundamental equations (1.4)--(1.7):
 $$\split(\widehat{\nabla}_XC)(Y,Z,W)&-(\widehat{\nabla}_YC)(X,Z,W)
=h(X,Z)h(SY,W)\\
 &-h(Y,Z)h(SX,W)+h(SY,Z)h(X,W)-h(SX,Z)h(Y,W)\endsplit\tag 3.1$$
 Evaluating this formula provides 
  $$\f_{ij}^k=0, \quad\text{if}\; \{i,j,k\}\neq\{1,2,3\},$$
  $$b-a=4\l(\f_{12}^3-\f_{21}^3),\tag 3.2$$
  $$c-a=4\l(\f_{13}^2-\f_{31}^2),\tag 3.3$$
 $$c-b=4\l(\f_{23}^1-\f_{32}^1),\tag 3.4$$
 $$e_1(\l)=e_2(\l)=e_3(\l)=0.\tag 3.5$$
 So adding (3.2) and (3.4) and substracting (3.3) yields
 $\f_{13}^2+\f_{32}^1-\f_{23}^1=0$. Let $\f_{13}^2=\a$,
 $\f_{32}^1=\b$, $\f_{21}^3=\g$, we therefore obtain $\a+\b+\g=0$, and 
 $$b-a=4\l\b,\quad a-c=4\l\g,\quad c-b=4\l\a;\tag 3.6$$
 $$\alignedat3 &\widehat{\nabla}_{e_1}e_1=0,&\qquad&\widehat{\nabla}_{e_1}e_2=-\a e_3
 ,&\qquad&\widehat{\nabla}_{e_1}e_3=\a e_2,\\
 &\widehat{\nabla}_{e_2}e_1=\g e_3,&\qquad&\widehat{\nabla}_{e_2}e_2=0,
 &\qquad&\widehat{\nabla}_{e_2}e_3=-\g e_1,\\
 &\widehat{\nabla}_{e_3}e_1=-\b e_2,&\qquad&\widehat{\nabla}_{e_3}e_2=\b e_1
 ,&\qquad&\widehat{\nabla}_{e_3}e_3=0.\endalignedat$$

From the fundamental equations (1.4)--(1.7) we also get the Gauss equations
 for $\widehat{R}$:
$$\widehat{R}(X,Y)Z=\frac{1}{2}\{h(Y,Z)SX-h(X,Z)SY+h(SY,Z)X-h(SX,Z)Y\}
-[K_X,K_Y]Z,
 \tag 3.7$$
 in which $\widehat{R}$ is the curvature tensor of $\widehat{\nabla}$.
 Note that
 $$[e_1,e_2]=\b e_3,\quad[e_3,e_1]=\g e_2,\quad[e_2,e_3]=\a
 e_1;$$
 $$[K_{e_1},K_{e_2}]=K_{e_1}K_{e_2}-K_{e_2}K_{e_1}=
 \pmatrix 0&-\l^2&0\\\l^2&0&0\\0&0&0\endpmatrix,$$
 $$[K_{e_1},K_{e_3}]=K_{e_1}K_{e_3}-K_{e_3}K_{e_1}=
 \pmatrix 0&0&-\l^2\\0&0&0\\\l^2&0&0\endpmatrix,$$
 $$[K_{e_2},K_{e_3}]=K_{e_2}K_{e_3}-K_{e_2}K_{e_3}=
 \pmatrix 0&0&0\\0&0&-\l^2\\0&\l^2&0\endpmatrix.$$
\smallskip
 Evaluating $\widehat{R}(e_1,e_2)e_1$ resp. $\widehat{R}(e_1,e_2)e_2$
 resp. $\widehat{R}(e_1,e_2)e_3$ yields
 $$e_1(\g)=0,\quad e_2(\a)=0,$$
$$\quad-\g\a-\b^2=\frac{1}{2}(a+b)+\l^2.\tag 3.8$$
 
Since $\{e_i\}$ and $\a$, $\b$, $\g$ admit cyclic symmetry,
 we get from (3.7) furthermore:
 $$e_1(\b)=e_3(\a)=e_2(\b)=e_3(\g)=0,$$
 $$-\a\b-\g^2=\frac{1}{2}(a+c)+\l^2,\tag 3.9$$
 $$-\b\g-\a^2=\frac{1}{2}(b+c)+\l^2.\tag 3.10$$

 \proclaim {Lemma 3.2} 
If $M^3$ admits a pointwise $Z_2\times Z_2$-symmetry, then
$$\a=b=c=0,\quad\l=\b=-\g=const.\neq0,\quad a=-4\l^2.$$
\endproclaim
 \demo{Proof}Since $\a+\b+\g=0$, we know that $\a$, $\b$, $\g$ are
 constant; (3.5) shows that $\l$ is constant, too. By (3.6),
 (3.8)--(3.10) we can express $a$, $b$, $c$ in terms of $\a$, $\b$,
 $\g$ and $\l$. Therefore, they are all constant.
From the difference of (3.8) and (3.10), resp. from (3.6), we get
 $$\b\g+\a^2-\g\a-\b^2=\frac{1}{2}(a-c)=2\l\g,\tag 3.11$$
 and similar from the difference of (3.10) and (3.9) we get
 $$\a\b+\g^2-\b\g-\a^2=\frac{1}{2}(b-a)=2\l\b.\tag 3.12$$ If we
 multiply (3.11) by $\b$ and (3.12) by $\g$ the resulting equations
 are equal, thus
 $$\aligned
 0&=2\a\b\g+\g^3+\b^3-\b\g^2-\a^2\g-\a^2\b-\b^2\g\\
 &=2\a\b\g+\g(\g^2-\a^2)+\b(\b^2-\a^2)-(\b\g^2+\b^2\g)\\
 &=2\a\b\g+\b\g(\a-\g)+\b\g(\a-\b)-\b\g(\b+\g)\\
 &=2\a\b\g+\b\g(\a+\a)-2\b\g(\b+\g)\\
 &=6\a\b\g.
 \endaligned$$
 If necessary, we can change the label of $\{e_i\}$, so that
 $\a=0$. From (3.6), $b=c$, and since $M^3$ is not a hypersphere, we
 have $\b=-\g\neq0$. Substituting into (3.11) we get $\l=\b$, then
 (3.10) yields $b=c=0$. Finally, $a=-4\l^2$ comes from (3.11).\qed
 \enddemo

 Therefore, we have the motion equations with respect to $\{e_i\}$,
 $$\alignat{3}
&D_{e_1} e_1 = &&&&\xi\\
&D_{e_1} e_2 =& && \l e_3& \\
&D_{e_1} e_3 =& &\l e_2&& \\
&D_{e_2} e_1 =&0 && &\\
&D_{e_3} e_1 =&0& & &\\
&D_{e_2} e_2 =&& & &\xi \\
&D_{e_2} e_3 =&2\l e_1& && \\
&D_{e_3} e_2 =& 2\l e_1&& & \\
&D_{e_3} e_3 =&& & &\xi\\
&D_{e_1} \xi=& 4\l^2 e_1&&&\\
&D_{e_2} \xi=&0&&&\\
&D_{e_3} \xi=&0&&&
\endalignat$$

Since $\widehat{\nabla}_{e_1}e_1=0$, locally we can choose a coordinate
$t$, such that $e_1=2\l\pt$. If we define two linear independent vector
fields
$$U=\frac{1}{2\l}e^{\frac{t}{2}}(e_2-e_3),\quad
V=\frac{1}{2\l}e^{-\frac{t}{2}}(e_2+e_3),$$ then it is easy to
verify that
$$[\,\pt,U\,]=[\,\pt,V\,]=[\,U,V\,]=0.$$
Thus there exist coordinates $\{u,v\}$, such that
$$U=\pu,\quad V=\pv,$$
and the motion equations in local coordinates $\{t,u,v\}$ are
$$\align
D_{\pt}\pt&=\frac{1}{4\l^2}\xi,\tag 3.13\\
D_{\pt}\pu=D_{\pt}\pv&=0,\tag 3.14\\
D_{\pu}\pu&=-2e^{t}\pt+\frac{e^{t}}{2\l^2}\xi,\tag 3.15\\
D_{\pv}\pv&=2 e^{t}\pt+\frac{e^{t}}{2\l^2}\xi,\tag 3.16\\
D_{\pu}\pv&=0,\tag 3.17\\
D_{\pt}\xi&=4\l^2\pt,\tag 3.18\\
D_{\pu}\xi=D_{\pv}\xi&=0.\tag 3.19\endalign$$

We can integrate these equations. If we denote the affine hypersurface
by $F(t,u,v):N\subset M^3\to\R^4$, then (3.14) implies that $F_t$ only
depends on $t$, so does $\xi$ by (3.19). By (3.18) and (3.13) we get
for $F_t$ that $(F_{t})_{tt}=(F_{tt})_{t}=F_t.$ So
$F_t=e^{t}\v_1+e^{-t}\v_2,$ in which $\v_1$, $\v_2$ are vectors in
$\R^4$, i.~e.
$$F=e^{t}\v_1-e^{-t}\v_2+\bold G(u,v),$$ 
and
$$\xi=4\l^2 e^{t}\v_1-4\l^2 e^{-t}\v_2.$$

From (3.15)-(3.17) we obtain the following differential equations for
$\bold G$: $\bold G_{uu}=-4\v_2$, $\bold G_{vv}=4\v_1$, $\bold
G_{uv}=0$, and the solution
$$\bold G(u,v)=2v^2\v_1-2u^2\v_2+v\v_3+u\v_4+\v_5.$$
 Hence,
 $$F=(e^{t}+2v^2)\v_1-(e^{-t}+2u^2)\v_2+v\v_3+u\v_4+\v_5.$$

 Note that $\v_1$, $\v_2$, $\v_3$ and $\v_4$ are linear independent
  since $\det (F_t, F_u, F_v, \xi) \neq 0$. By an affine
  transformation we can translate $\v_5$ to zero and map $\v_1$,
  $\v_2$, $\v_3$, $\v_4$ to $(1,0,0,0)$, $(0,-1,0,0)$, ($0,0,2,0)$,
  $(0,0,0,2)$, resp., and obtain

$$F(t,u,v)=\left(e^{t}+2v^2,e^{-t}+2u^2,2v,2u\right).$$

Thus we obtain
 \proclaim{Theorem 3.1}  $F:M^3\to\R^4$ is a locally strictly
 convex affine hypersurface admitting pointwise $Z_2\times Z_2$-symmetry,
 if and only if it is affine equivalent to
$$
(x_1-\frac12 x_3^2)(x_2-\frac12 x_4^2)=1.
$$
 \endproclaim
\remark{Remark} 
The affine hypersurface above is one of the two types of locally
homogeneous affine hypersurfaces with rank one shape operator,
classified in \cite{3}.
\endremark

 \head 4. Hypersurfaces admitting $SO(2)$- or $Z_3$-symmetry\endhead

 What is left, is the classification of hypersurfaces admitting pointwise 
$SO(2)$- or $Z_3$-symmetry. This will be done as follows. Let $M^3$ be an 
affine hypersurface admitting pointwise 
$SO(2)$- or $Z_3$-symmetry. First, we will obtain information
about the coefficients of the connection from the basic equations
of Gauss, Codazzi and Ricci. In particular, it follows that such a
hypersurface $M^3$ admits a warped product structure. Then following
essentially the same approach as in \cite {12}, we classify such
hypersurfaces by showing how they can be constructed starting from
$2$-dimensional positive definite affine spheres. 

Affine hyperspheres were classified in \cite {12}, thus we could exclude them 
in the following. Anyhow, it is not much harder to handle the general case.

\bigskip
 \subhead 4.1. Structure equations and integrability
 conditions\endsubhead
 \subsubhead 4.1.1 An adapted frame \endsubsubhead
 \bigskip

 We recall from Corollary 2.3, Case (iii) and (vii): At every point $p$ of 
$M^3$ there exists
a basis $\{e_1,e_2,e_3\}$ which is orthonormal with respect to the
affine metric $h$ such that the difference tensor $K$ and the
shape operator $S$ are respectively given by:
 $$\aligned
  K_{e_1}&=\pmatrix 2\l &0&0\\ 0& -\l & 0 \\ 0&0& -\l
  \endpmatrix, \quad
  K_{e_2}=\pmatrix 0 & -\l&0\\ -\l & \mu & 0 \\ 0&0& -\mu
  \endpmatrix, \\
K_{e_3}&=\pmatrix 0 & 0 & -\l\\ 0 & 0 & -\mu \\ -\l &-\mu& 0
  \endpmatrix,\quad
S=\pmatrix a& 0 &0 \\ 0& b &0 \\ 0& 0 & b
  \endpmatrix.
\endaligned\tag 4.1$$
 We have that $\l>0$. Moreover, $\mu$ vanishes if and only
if the symmetry group is $SO(2)$, i.e.
the form of $K$ and $S$ remains invariant under rotations in the
$e_2e_3$-plane.  In case that $\mu$ is different from zero, the
group $Z_3$ of rotations leaving $K$ and $S$ invariant is
generated by $R_1$. Moreover, if $\mu=\sqrt{2} \l$, then $a\neq b$. 
We want to extend this basis differentiably.

We define the Ricci tensor of the connection $\widehat \nabla$ by:
$$\ricm(X,Y)=\trace\{ Z \mapsto \widehat R(Z,X)Y\}.$$
It is well known that $\ricm$ is a symmetric operator. Then, we
have
 \proclaim{Lemma 4.1}
 Let $p \in M^3$ and $\{e_1,e_2,e_3\}$ the basis constructed
  earlier. Then
  $$\alignat{2}
  &\ricm(e_1,e_1) = (a+b)+ 6 \l^2, \qquad\quad &&\ricm(e_1,e_2)=0, \\
  &\ricm(e_3,e_1)=0,\qquad\quad &&
\ricm(e_2,e_2)=\tfrac{1}{2}a +\tfrac{3}{2} b +2(\l^2+\mu^2),\\
  &\ricm(e_2,e_3)=0,&&\ricm(e_3,e_3)=\tfrac{1}{2}a +\tfrac{3}{2}
  b +2(\l^2+\mu^2).
\endalignat$$
\endproclaim
 \demo{Proof}We use the Gauss equation (3.6) for $\widehat R$.
It follows that
 $$\align
  \widehat R(e_2,e_1)e_1&=\tfrac{1}{2}(a+b) e_2-K_{e_2}(2\l
  e_1)+K_{e_1}(-\l e_2) \\
  &=(\tfrac{1}{2}(a+b)+3 \l^2) e_2,\\
  \widehat R(e_3,e_1)e_1&=\tfrac{1}{2}(a+b) e_3-K_{e_3}(2\l e_1)
  +K_{e_1}(-\l e_3) \\
  &=(\tfrac{1}{2}(a+b)+3 \l^2) e_3,\\
  \widehat R(e_3,e_1)e_2&=-K_{e_3}(-\l e_2)+K_{e_1}(-\mu e_3)=0.
\endalign$$
From this it immediately follows that
$$\ricm(e_1,e_1) = (a+b)+6 \l^2$$
and
$$\ricm(e_1,e_2)=0.$$
The other equations follow by similar computations.\qed
\enddemo
 Now, we want to show that the basis we have constructed at each
point $p$ can be extended differentiably to a neighborhood of the
point $p$ such that at every point the components of $S$ and $K$
with respect to the frame $\{e_1,e_2,e_3\}$ have the previously
described form.
 \proclaim{Lemma 4.2}Let $M^3$ be an affine hypersurface of $\R^4$ which
  admits a pointwise $Z_3$-symmetry or a pointwise
  $SO(2)$-symmetry. Let $p \in M$. Then there exists a frame
  $\{e_1,e_2,e_3\}$ defined in a neighborhood of the point $p$ such that
  the components of $K$ and $S$ are respectively given by:
  $$\aligned
  K_{e_1}&=\pmatrix 2\l &0&0\\ 0& -\l & 0 \\ 0&0& -\l
  \endpmatrix, \quad
  K_{e_2}=\pmatrix 0 & -\l&0\\ -\l & \mu & 0 \\ 0&0& -\mu
  \endpmatrix, \\
K_{e_3}&=\pmatrix 0 & 0 & -\l\\ 0 & 0 & -\mu \\ -\l &-\mu& 0
  \endpmatrix,\quad
S=\pmatrix a& 0 &0 \\ 0& b &0 \\ 0& 0 & b
  \endpmatrix.
\endaligned$$
 \endproclaim
 \demo{Proof}First we want to show that at every point the vector
  $e_1$ is uniquely defined and differentiable. We introduce a
  symmetric operator $\hat A$ by:
  $$\ricm(Y,Z)= h(\hat A Y,Z).$$
  Clearly $\hat A$ is a differentiable operator on $M$. On the set
of points where $\tfrac{1}{2}(a-b) +4\l^2 -2 \mu^2\neq 0$, the
operator has two distinct eigenvalues. The eigendirection which
corresponds with the $1$-dimensional eigenvalue corresponds with
the vector field $e_1$.

On the set of points where $\tfrac{1}{2}(a-b) +4\l^2 -2 \mu^2= 0$,
we can assume that $a \neq b$. (If $a=b$ then $\mu \neq \sqrt2 \l$ or 
$\mu=0$ by Corollary 2.3.) In
this case, the differentiable operator $S$ has two distinct
eigenvalues (cf. (4.1)) and $e_1$ is uniquely determined as the
eigendirection corresponding to the $1$-dimensional eigenvalue.
This shows that taking at every point $p$ the vector $e_1$ yields
a differentiable vector field.

To show that $e_2$ and $e_3$ can be extended differentiably, we
consider two cases. First we assume that $M$ admits a pointwise
$SO(2)$-symmetry. In that case we have that $\mu=0$ and we take
for $e_2$ and $e_3$ arbitrary orthonormal differentiable local
vector fields which are orthogonal to the vector field $e_1$. In
case that $M$ admits a pointwise $Z_3$-symmetry we proceed as
follows. We start by taking arbitrary orthonormal differentiable
local vector fields $u_2$ and $u_3$ which are orthogonal to the
vector field $e_1$. It is then straightforward to check that we
can write
$$
\alignat{2}
  &K_{u_2}u_2= -\l e_1 &+&\nu_1 u_2 +\nu_2 u_3\\
  &K_{u_2}u_3= &&\nu_2 u_2 -\nu_1 u_3\\
  &K_{u_3}u_3= -\l e_1 &-&\nu_1 u_2 -\nu_2 u_3\endalignat
$$
differentiable functions $\nu_1$ and $\nu_2$ with $\nu_1^2
+\nu_2^2\neq 0$. Therefore, if necessary by interchanging the role
of $u_2$ and $u_3$, we may assume that in a neighborhood of the
point $p$, $\nu_1\neq 0$.  Rotating now over an angle $\theta$,
thus defining $$\align
&e_2 =\cos \theta u_2 +\sin \theta u_3,\\
&e_3 =-\sin \theta u_2 + \cos \theta u_3,
\endalign$$
we get that $$\align
  h(K(e_2,e_2),e_2)&=(\cos^3 \theta -3 \cos \theta \sin^2 \theta)
  \nu_1 + (-\sin^3 \theta +3 \cos^2 \theta \sin \theta)\nu_2\\
  &=\cos 3 \theta \nu_1 +\sin 3 \theta \nu_2\\
  h(K(e_3,e_3),e_3)&=(-\sin^3 \theta +3 \cos^2 \theta \sin \theta)
  \nu_1 + (-\cos^3 \theta +3 \cos \theta \sin^2 \theta)\nu_2\\
  &=\sin 3 \theta \nu_1 - \cos 3 \theta \nu_2.
\endalign$$
Therefore, taking into account the symmetries of $K$, in order to
obtain the desired frame, it is sufficient to choose $\theta$ in
such a way that
$$
\sin 3 \theta \nu_1 - \cos 3 \theta \nu_2=0,
$$
and $\cos 3 \theta \nu_1 +\sin 3 \theta \nu_2>0$. As this is always
possible, the proof is completed.\qed 
\enddemo 
\remark{Remark} 
It actually follows from the proof of the previous lemma that the
vector field $e_1$ is globally defined on $M^3$, and therefore the
function $\l$, too. This in turn implies that the functions $\mu$ (as
it can be expressed in terms of $\l$ and the Pick invariant $J$
(cf. \cite{9}, Prop. 9.3) and $J$ is either identically zero or
nowhere zero), $a$ and $b$ (as it can be expressed in terms of the
mean curvature and $a$, cf. \cite{9}, Def. 3.4) are globally defined
functions on the affine hypersurface $M^3$.
\endremark
\smallskip
From now on we always will work with the local frame constructed
in the previous lemma. We introduce the connection coefficients with
respect to this frame by $\widehat \nabla_{e_i} e_j =\sum_{k=1}^3
\f_{ij}^{k} e_k$, as we did before, and have the usual symmetries.
\bigskip
 \subsubhead 4.1.2 Codazzi equations for $K$\endsubsubhead
\bigskip
An evaluation of the Codazzi equations for $K$ (cp.(3.1)) using the
CAS mathematica (cf. \cite{10}) results in the following equations:
 $$\allowdisplaybreaks\align
e_2(\l)=2 \l \f_{11}^2, &\qquad (eq. 1)\tag 4.2\\
\mu \f_{11}^3 =4 \l \f_{21}^3, & \qquad(eq. 1)\tag 4.3\\
e_1(\l)= \tfrac{1}{2} (b-a)-\mu \f_{11}^2- 4\l \f_{21}^2,
&\qquad(eq. 1)\tag 4.4\\
e_1(\l)=\tfrac{1}{2} (b-a)+\mu \f_{11}^2 -4\l \f_{31}^3, &
\qquad(eq. 2)\tag 4.5
\\
e_3(\l)=2 \l \f_{11}^3, & \qquad(eq. 2)\tag 4.6\\
\mu \f_{11}^3=4 \l \f_{31}^2, & \qquad(eq. 2)\tag 4.7\\
e_1(\mu)+e_2(\l)= 3\l \f_{11}^2 -\mu \f_{21}^2, & \qquad(eq. 3)\tag 4.8\\
0=-\l \f_{11}^3 + 3\mu \f_{12}^3 -\mu \f_{21}^3, & \qquad(eq. 3)\tag 4.9\\
e_3(\l)= -\mu (\f_{21}^3 +\f_{31}^2), & \qquad(eq. 4)\tag 4.10\\
e_3(\mu)=3\mu \f_{22}^3 -\l (\f_{21}^3-3\f_{31}^2), &\qquad (eq. 4)\tag 4.11\\
e_2(\mu)=-\l (\f_{21}^2-\f_{31}^3) -3 \mu \f_{32}^3, & \qquad(eq. 4)\tag 4.12\\
e_1(\mu)=-\l\f_{11}^2 -\mu \f_{31}^3, &\qquad(eq. 5)\tag 4.13\\
e_3(\l)=3 \l \f_{11}^3 +\mu(3\f_{12}^3+\f_{31}^2), &\qquad(eq. 5)\tag 4.14\\
e_2(\l)=\mu (\f_{21}^2-\f_{31}^3), &\qquad(eq. 6)\tag 4.15\\
e_3(\mu)=3\mu \f_{22}^3 +\l (3 \f_{21}^3-\f_{31}^2), &\qquad (eq. 6)\tag 4.16\\
e_2(\l)-e_1(\mu)=\l \f_{11}^2 +\mu \f_{21}^2, &\qquad(eq. 7)\tag 4.17\\
4 \l (\f_{21}^3 -\f_{31}^2) =0,&\qquad(eq. 8)\tag 4.18\\
e_3(\l)=\l \f_{11}^3- \mu (3\f_{12}^3+\f_{31}^2). &\qquad(eq.
9)\tag 4.19
\endalign$$

In the above expressions, the equation numbers refer to
corresponding equations in the mathematica program.  In order to
simplify the above equations, we now distinct two cases.
 \proclaim{Lemma 4.3}
 An evaluatin of the Codazzi equations for
  $K$ gives:
$$ \f_{11}^2=0,\quad \f_{11}^3=0,\quad \f_{21}^3=0,\quad \f_{31}^2=0,
  \quad \f_{21}^2= \f_{31}^3 =:\eta,$$
  $$e_1(\l)=\frac{1}{2}(b-a) -4 \l \eta,\quad e_2(\l)=0=e_3(\l).
$$
If $\mu\neq 0$, we get in addition that $\f_{12}^3=0$ and
$$
e_1(\mu)= -\mu \eta, \quad e_2(\mu)= -3 \mu \f_{32}^3,\quad
e_3(\mu)=3\mu \f_{22}^3.\tag 4.20
$$
 \endproclaim
 \demo{Proof}
 First, we assume that $\mu=0$ (thus $\l\neq 0$). In that
  case, it follows from (4.3) (resp. (4.7)) that
  $\f_{21}^3=0$ (resp.  $\f_{31}^2=0$), whereas (4.12) implies that
  $\f_{21}^2 =\f_{31}^3$. As it now follows from (4.10) and
  (4.15) that $e_2(\l)=e_3(\l)=0$, (4.2) and
  (4.6) imply that $\f_{11}^2=\f_{11}^3=0$. Finally (4.5)
  now reduces to
$$e_1(\l)= \tfrac{1}{2} (b-a)- 4\l \f_{21}^2.$$

Next, we want to deal with the case that $\mu \neq 0$. First it
follows from (4.2), (4.15) and (4.8), taking also into account
(4.13), that
$$e_2(\l)=2 \l \f_{11}^2=\mu (\f_{21}^2-\f_{31}^3)=4 \l \f_{11}^2
-\mu (\f_{21}^2-\f_{31}^3).$$ Therefore we get by (4.4) and (4.5)
that $\f_{21}^2=\f_{31}^3$ and thus $e_2(\l)=0=\f_{11}^2$. From
(4.10) and (4.6) it follows that
$$e_3(\l)=-\mu (\f_{21}^3 +\f_{31}^2)=2 \l \f_{11}^3.$$
From (4.18), (4.3) and the previous equation it follows that
$\f_{21}^3=\f_{31}^2=\f_{11}^3=0$ and $e_3(\l)=0$. From (4.9) it
then follows that $\f_{12}^3=0$.  From (4.8), (4.12) and (4.11) we
obtain the equations for $e_i(\mu)$, $i=1,2,3$ and from (4.4) it
follows that
$$e_1(\l)= \tfrac{1}{2} (b-a)- 4\l \f_{21}^2. \qed$$\enddemo
As a direct consequence we write down the Levi-Civita connection:
\proclaim{Lemma 4.4} $$\align
   \widehat{\nabla}_{e_1}e_1 &= 0,\\
   \widehat{\nabla}_{e_1}e_2 &= \f_{12}^3 e_3,\\
   \widehat{\nabla}_{e_1}e_3 &= -\f_{12}^3 e_2,\\
   \widehat{\nabla}_{e_2}e_1 &= \eta e_2,\\
   \widehat{\nabla}_{e_2}e_2 &= -\eta e_1 + \f_{22}^3 e_3,\\
   \widehat{\nabla}_{e_2}e_3 &= -\f_{22}^3 e_2,\\
   \widehat{\nabla}_{e_3}e_1 &= \eta e_3,\\
   \widehat{\nabla}_{e_3}e_2 &= \f_{32}^3 e_3,\\
   \widehat{\nabla}_{e_3}e_3 &= -\eta e_1- \f_{32}^3 e_2,
   \endalign$$
   where in case that $\mu \neq 0$, we have in addition that
$\f_{12}^3=0$.
\endproclaim
\bigskip
\subsubhead 4.1.3 Gauss for $\nabla$\endsubsubhead
\bigskip
Taking into account the previous results, we then proceed with an
evaluation of the Gauss equations (1.4) for $\nabla$:
$$
\nabla_X \nabla_Y Z- \nabla_Y \nabla_X Z - \nabla_{[X,Y]} Z =
h(Y,Z)SX - h(X,Z)SY,
$$
again using the CAS mathematica (cf. \cite{10}).
This results amongst others in the following equations (cf.
equations 11, 13, 14 and 16 in the mathematica program):
$$\align
e_1(\eta)& = -\eta^2 -3 \l^2 -\frac{1}{2}(a+b), \\
e_2(\eta)& = 0,\\
e_3(\eta)& =0.
\endalign$$
\bigskip
\subsubhead 4.1.4 Codazzi for $S$\endsubsubhead
\bigskip
An evaluation of the Codazzi equations (1.6) for $S$:
$$
(\nabla_X S)(Y) = (\nabla_Y S)(X)
$$
by mathematica (cf. \cite{10}, equations 20 - 22)) then yields:
$$\align
  e_1(b) &= (\l-\eta)(b-a),\\
e_2(b) &= 0, \\
e_3(b) &= 0, \\
e_2(a) &= 0, \\
e_3(a) &= 0.
\endalign$$
\bigskip
\subsubhead 4.1.5 Structure equations\endsubsubhead
\bigskip
Summarized we have obtained the structure equations (cf. (1.1),
(1.2) and (1.3)):
$$\alignat{4}
&D_{e_1} e_1 = &2\l e_1&&&+ \xi,\tag4.21 \\
&D_{e_1} e_2 =& &-\l e_2& + \f_{12}^3 e_3,&\tag4.22 \\
&D_{e_1} e_3 =& &-\f_{12}^3 e_2&-\l e_3,& \tag4.23\\
&D_{e_2} e_1 =& &(\eta -\l) e_2,& & \tag4.24\\
&D_{e_3} e_1 =&& &(\eta -\l) e_3,  &\tag4.25\\
&D_{e_2} e_2 =&-(\eta+\l) e_1& +\mu e_2& + \f_{22}^3 e_3 &+\xi, \tag4.26\\
&D_{e_2} e_3 =&&-\f_{22}^3 e_2 &-\mu e_3,& \tag4.27\\
&D_{e_3} e_2 =& &&(\f_{32}^3-\mu) e_3, & \tag4.28\\
&D_{e_3} e_3 =&-(\eta+\l) e_1& -(\f_{32}^3-\mu) e_2&
&+\xi,\tag4.29
\endalignat$$
$$\alignat{4}
&D_{e_1} \xi=& -a e_1,&&&\tag4.30\\
&D_{e_2} \xi=&& -b e_2,&&\tag4.31\\
&D_{e_3} \xi=&&& -b e_3,&\tag4.32
\endalignat$$
Moreover, the functions $a$, $b$, $\l$ and $\eta$ are all constant
in the $e_2$- and $e_3$-directions and the $e_1$-derivatives are
determined by (cf. Subsection 4.1.3 and 4.1.2 and Lemma 4.3):
$$\align
&e_1(b) = (\l-\eta)(b-a),\tag4.33\\
&e_1(\eta) = -\eta^2 -3 \l^2 -\tfrac{1}{2}(a+b), \tag4.34\\
&e_1(\l)= - 4\l \eta-\tfrac{1}{2} (a-b).\tag4.35
\endalign$$
 \bigskip
\subhead 4.2. The classification \endsubhead
 \bigskip
 As the vector field $e_1$ is globally defined, we can define the
distributions $H_1=\operatorname{Span}\{e_1\}$ and
$H_2=\operatorname{Span}\{e_2,e_3\}$. In the next lemmas we will
investigate some properties of these distributions following from
Lemma 4.4. For the terminology we refer to \cite{8}.
 \proclaim{Lemma 4.5}
 The distribution $H_1$ is autoparallel with respect to
$\widehat\nabla$.
 \endproclaim
 \demo{Proof}From $\widehat{\nabla}_{e_1} e_1 =0$ the claim follows
  immediately.
 \qed\enddemo
 \proclaim{Lemma 4.6}
 The distribution $H_2$ is spherical with mean curvature normal
  $H=-\eta e_1$.
 \endproclaim
 \demo{Proof}For $H=-\eta e_1\in H_1=H_2^{\perp}$ we have
  $h(\widehat{\nabla}_{e_a} e_b, e_1)= h(e_a, e_b) h(H,e_1)$ for all $a,b
  \in \{2,3\}$, and $h(\widehat{\nabla}_{e_a} H, e_1)= h(-e_a(\eta) e_1
  - \eta \widehat{\nabla}_{e_a} e_1, e_1)=0$.
 \qed\enddemo
 \remark{Remark}
 $\eta\; (=\f_{21}^2=\f_{31}^3)$ is independent of the particular 
  choice of the orthonormal basis $\{e_2, e_3\}$. It therefore is a 
  globally defined function on $M^3$.
 \endremark
\smallskip
 We introduce a coordinate function $t$ by $\pt:=e_1$. Using the
previous lemma, according to \cite{5}, we get:
 \proclaim{Lemma 4.7}
 $(M,h)$ admits a warped product structure
  $M^3=\R\times_{e^f}N^2$ with $f:\R \to \R$
  satisfying
$$
\frac{\partial f}{\partial t}=\eta.\tag 4.36
$$
 \endproclaim
 \remark{Remark}
 $b$, $\eta$ and $\l$ are functions of $t$, they satisfy
  by (4.33), (4.34) and (4.35):
 $$\align
   \frac{\partial b}{\partial t}&= (\l-\eta)(b-a),\\
   \frac{\partial \eta}{\partial t}&=-\eta^2-3\l^2-\tfrac 12 (a+b),\\
   \frac{\partial \l}{\partial t}&=-4 \eta \l +\tfrac{1}{2}
   (b-a).
 \endalign$$
 \endremark
 To compute the curvature of $N^2$ we use the gauss equation
(3.6) and obtain:
$$
   K(N^2)=e^{2f}(b-\l^2+2\mu^2+\eta^2),\tag 4.37
$$
which we verify by a straightforward computation is indeed
independent of $t$.

Our first goal is to find out how $N^2$ is immersed in $\R^4$,
i.~e. to find an immersion independent of $t$. A look at the
structure equations (4.21) - (4.32) suggests to start with a
linear combination of $e_1$ and $\xi$.

We will solve the problem in two steps. First we define $v:=A e_1
+\xi$ for some function $A$ on $M^3$. Then $\pt v=c v$ iff $c=A$
and $\pt A= A^2 -2\l A+ a$, and $A:=-(\eta+\l)$ solves the latter
differential equation. Next we define a positive function $\tau$
 on $\R$ as solution of the differential equation:
$$\tfrac{\partial}{\partial t} \tau = \tau(\eta+\l) \tag 4.38$$
 with initial condition $\tau(t_0)>0$. Then $\pt(\tau v)=0$ and by (4.24),
(4.31), (4.25) and (4.32) we get (since $\tau$, $\eta$ and $\l$ only
depend on $t$): $$\align
D_{e_1}(\tau(-(\eta+\l) e_1 +\xi))&=0,\tag 4.39\\
D_{e_2}(\tau(-(\eta+\l) e_1 +\xi))&=-\tau(b +\eta^2 -\l^2)
e_2,\tag4.40\\
D_{e_3}(\tau(-(\eta+\l)e_1 +\xi))&=-\tau(b +\eta^2 -\l^2) e_3.\tag4.41
\endalign$$
\proclaim{Lemma 4.8}
 Define $\nu:=b +\eta^2 -\l^2$ on $\R$. $\nu$ is globally defined,
  $\pt(e^{2f} \nu)=0$ and $\nu$ vanishes identically or nowhere on $\R$.
 \endproclaim
 \demo{Proof}Since $0=\pt K(N^2) = \pt(e^{2f}(\nu+2\mu^2))$
  (cf. (4.20) and (4.36)) and $\pt(e^{2f}2 \mu^2)=0$, we
  get that $\pt(e^{2f} \nu)=0$. Thus $\pt\nu=-2 (\pt f)\nu= -2\eta\nu$.
 \qed\enddemo

 Now we consider different cases depending on the behaviour of
$\nu$.

\specialhead 4.2.1 The first case: $\nu \neq 0$ on
$M^3$\endspecialhead

We may, by translating $f$, i.e. by replacing $N^2$ with a
homothetic copy of itself, assume that $e^{2f} \nu =\epsilon$,
where $\epsilon =\pm 1$.

\proclaim{Lemma 4.9} $\varPhi:=\tau (-(\eta +\l) e_1 +\xi)\colon
M^3 \to \R^4$ induces a proper affine sphere structure, say
$\tilde{\phi}$, mapping $N^2$ into a 3-dimensional linear subspace
of $\R^4$. $\tilde{\phi}$ is part of a quadric iff $\mu =0$.
\endproclaim
\demo{Proof} By (4.40) and (4.41) we have $\varPhi_*(e_i)= -\tau
\nu e_i$ for $i= 2,3$. A further differentiation, using (4.26)
($\tau$ and $\nu$ only depend on $t$), gives:
$$\align
D_{e_2} \varPhi_*(e_2)& = -\tau (b +\eta^2-\l^2) D_{e_2}e_2\\
&= -\tau (b +\eta^2-\l^2) (-(\eta+\l) e_1 +\mu e_2 + \f_{22}^3 e_3 +\xi)\\
&=\mu\varPhi_*(e_2)+\f_{22}^3 \varPhi_*(e_3) -(b+\eta^2 -\l^2) \varPhi\\
&=\mu\varPhi_*(e_2)+\f_{22}^3 \varPhi_*(e_3) -\epsilon e^{-2f}
\varPhi.
\endalign$$
Similarly, we obtain the other derivatives,using (4.27)-(4.29),
thus:
 $$\alignedat{3}
 D_{e_2} \varPhi_*(e_2)&= &\mu \varPhi_*(e_2) &+
\f_{22}^3 \varPhi_*(e_3)&- e^{-2f}\epsilon \varPhi, \\ 
D_{e_2} \varPhi_*(e_3)&=& -\f_{22}^3 \varPhi_*(e_2) &-\mu
\varPhi_*(e_3),& \\ 
D_{e_3} \varPhi_*(e_2)&=& &(\f_{32}^3-\mu)
\varPhi_*(e_3),& \\ D_{e_3} \varPhi_*(e_3)&=&
-(\f_{32}^3+\mu)\varPhi_*(e_2) &&- e^{-2f}\epsilon \varPhi, \\ 
D_{e_i}\varPhi &=& -\tau e^{-2f}\epsilon e_i, & &\\
&=&\varPhi_*(e_i),&\qquad i=2,3.&
\endalignedat$$

The foliation at $f=f_0$ gives an immersion of $N^2$ to $M^3$, say
$\pi_{f_0}$. Therefore, we can define an immersion of $N^2$ to
$\R^4$ by $\tilde{\phi}:=\varPhi\circ\pi_{f_0}$, whose motion
equations are exactly the equations above when $f=f_0$. Hence, we
know actually $\tilde{\phi}$ maps $N^2$ into
$\sp\{\varPhi_*(e_2),\varPhi_*(e_3),\varPhi\}$, an affine
hyperplane of $\R^4$. $\pt\varPhi=0$ implies
$\varPhi(t,u,v)=\tilde{\phi}(u,v)$.

 We can read off
the coefficients of the difference tensor $K^{\tilde{\phi}}$ of
$\tilde{\phi}$ (cf. (1.1) and (1.3)):
$(K^{\tilde{\phi}})_{22}^2=\mu$, $(K^{\tilde{\phi}})_{23}^3=-\mu$,
$(K^{\tilde{\phi}})_{22}^3=0= (K^{\tilde{\phi}})_{33}^3$, and see
that $\trace (K^{\tilde{\phi}})_X$ vanishes. The affine metric
introduced by this immersion corresponds with the metric on $N^2$.
Thus $-\epsilon \varPhi$
 is the affine normal of $\tilde{\phi}$ and $\tilde{\phi}$ is a
 proper affine sphere with mean curvature $\epsilon$.
Finally the vanishing of the difference tensor characterizes
quadrics.\qed
\enddemo
Our next goal is to find another linear combination of $e_1$ and
$\xi$, this time only depending on $t$. (Then we can express $e_1$
in terms of $\varPhi$ and some function of $t$.)
 \proclaim{Lemma 4.10}
  Define $\delta := b e_1 +(\eta -\l) \xi$. Then there exist a constant vector
  $C \in \R^4$ and a function $g(t)$ such that
$$ \delta(t)= g(t) C.$$
\endproclaim
 \demo{Proof}
 Using (4.24) and (4.31) resp. (4.25) and
  (4.32) we obtain that $D_{e_2}\delta = 0=D_{e_3} \delta$. Hence $\delta$
  depends only on the variable $t$. Moreover, we get by
  (4.33), (4.21), (4.34), (4.35) and (4.30)
  that
 $$\align
  \pt\delta&=D_{e_1} ( b e_1+(\eta-\l)\xi)\\
  &=(\l-\eta)(b-a) e_1 + 2 b \l e_1+b \xi -(\eta-\l)a e_1\\
  &\qquad +(-\eta^2 -3 \l^2 -b+ 4\l \eta)\xi\\
  &=(3 \l-\eta)(b e_1 +(\eta -\l)\xi)\\
  &=(3 \l -\eta) \delta.
\endalign$$
This implies that there exists a constant vector $C$ in $\R^4$ and
a function $g(t)$ such that $\delta(t)=g(t)C$.\qed
 \enddemo
  Combining $\varPhi$ and $\delta$ we obtain for $e_1$ (cf.
Lemma 4.9 and 4.10) that
$$\align
e_1(t,u,v)&= \tfrac{1}{\tau\nu}(\tau gC
-(\eta-\l)\varPhi(t,u,v))\\
&=\tfrac{g}{\nu}C
-\tfrac{1}{\tau\nu}(\eta-\l)\tilde{\phi}(u,v).\tag4.42
\endalign
$$

In the following we will use for the partial derivatives the
abbreviation $F_x:= \frac{\partial}{\partial x}F $, $x=t,u,v$.
\proclaim{Lemma 4.11}
  $$\align&F_t = \frac{g}{\nu}C -\pt(\frac{1}{\tau \nu})\tilde{\phi},\\
&F_u = -\frac{1}{\tau\nu} \tilde{\phi}_u,\\
&F_v = -\frac{1}{\tau\nu} \tilde{\phi}_v.\endalign$$

\endproclaim
 \demo{Proof}
 As by (4.38) and Lemma 4.8 $\pt
  \frac{1}{\tau\nu}= \frac{1}{\tau\nu}(\eta -\l)$, we obtain the equation for $F_t
  =e_1$ by (4.42). The other equations follow from (4.40)
  and (4.41).\qed
 \enddemo
 It follows by the uniqueness theorem of first order differential
equations and applying a translation that we can write
$$F(t,u,v)= \tilde{g}(t) C -\frac{1}{\tau\nu}(t) \tilde{\phi}(u,v)$$
for a suitable function $\tilde{g}$ depending only on the variable
$t$. Since $C$ is transversal to the image of $\tilde{\phi}$ (cf.
Lemma 4.9 and 4.10), we obtain that after applying an equiaffine
transformation we can write: $F(t,u,v) =(\gamma_1(t), \gamma_2(t)
\phi(u,v))$, in which $\tilde{\phi}(u,v)=(0,\phi(u,v))$. Thus we
have proven the following:
 \proclaim{Theorem 4.1}
  Let $M^3$ be a positive definite affine hypersurface of $\R^4$ which
  admits a pointwise $SO(2)$- or $Z_3$-symmetry and with the
  globally defined function $(b +\eta^2 -\l^2)$ not identically zero
  on $M^3$. Then $M^3$ is affine equivalent to
  $$F:I\times N^2\to \R^4:(t,u,v)\mapsto (\gamma_1(t),
  \gamma_2(t) \phi(u,v)),$$
  where
  $\phi: N^2 \to \R^3$ is an elliptic or hyperbolic affine sphere
  and $\gamma:I\to \R^2$ is a curve.\hfill \newline
  Moreover,
  if $M^3$ admits a pointwise $SO(2)$-symmetry then $N^2$ is either an
  ellipsoid or a hyperboloid.
\endproclaim

In the next theorem we deal with the converse.  
\proclaim{Theorem 4.2}
 Let $\phi:N^2 \to \R^3$ be an elliptic or hyperbolic affine sphere
 \rom(scaled such that the absolute value of the mean curvature equals
 $1$\rom ). Let $\gamma: I \to \R^2$ be a curve such that
  $$F(t,u,v)=(\gamma_1(t), \gamma_2(t) \phi(u,v)),$$ defines a
  positive definite affine hypersurface. Then $F$ admits a pointwise
  $Z_3$- or $SO(2)$-symmetry.
\endproclaim
 \demo{Proof}
 We have
 $$\align
  &F_t=(\gamma_1',\gamma_2' \phi),\\
  &F_u=(0,\gamma_2 \phi_u),\\
  &F_v=(0,\gamma_2 \phi_v),\\
  &F_{tt}=(\gamma_1'',\gamma_2'' \phi)=\tfrac{(\gamma_2''\gamma_1'-
  \gamma_1''\gamma_2')}{\gamma_1'}(0,\phi)+\tfrac{\gamma_1''}{\gamma_1'}
  F_t,\\
  &F_{ut}=\tfrac{\gamma_2'}{\gamma_2} F_u,\\
  &F_{vt}=\tfrac{\gamma_2'}{\gamma_2} F_v,\\
  &F_{uu}=(0,\gamma_1\phi_{uu}),\\
  &F_{uv}=(0,\gamma_1\phi_{uv}),\\
  &F_{vv}=(0,\gamma_1\phi_{vv}).
\endalign$$
This implies that $F$ defines a nondegenerate affine immersion
provided that $\gamma_2\gamma_1
\gamma_1'(\gamma_2''\gamma_1'-\gamma_1'' \gamma_2') \neq 0$.  We
moreover see that this immersion is definite provided that the
affine sphere is hyperbolic and $\gamma_1
\gamma_1'(\gamma_2''\gamma_1'-\gamma_1'' \gamma_2')>0$ or when the
proper affine sphere is elliptic and $\gamma_1
\gamma_1'(\gamma_2''\gamma_1'-\gamma_1'' \gamma_2')<0$.  As the
proof in both cases is similar, we will only treat the first case
here. An evaluation of the conditions for the affine normal $\xi$
($\xi_t$, $\xi_u$ , $\xi_v$ are tangential and
$\det(F_t,F_u,F_v,\xi)=\sqrt{\det
  h}$) leads to:
$$
\xi=\alpha(t) (0,\phi(u,v)) +\beta(t) F_t,
$$
where $(\gamma_2''\gamma_1'-\gamma_1'' \gamma_2')\gamma_1^2 =
\gamma_2^4 (\gamma_1')^3 \alpha^5$ and $\alpha'
+\tfrac{\beta(\gamma_2''\gamma_1'-\gamma_1''\gamma_2')}{\gamma_1'}=0$.
Taking $e_1$ in the direction of $F_t$, we see that $F_u$ and
$F_v$ are orthogonal to $e_1$. It is also clear that $S$
restricted to the space spanned by $F_u$ and $F_v$ is a multiple
of the identity, and $S(F_t)=a F_t$, since $S$ is symmetric.
Moreover, we have that
 $$\align
  (\nabla h)(F_t,F_u,F_u)
  &=(\tfrac{\gamma_1'}{\gamma_1}-\tfrac{\alpha'}{\alpha}-2
  \tfrac{\gamma_2'}{\gamma_2})h(F_u,F_u),\\
  (\nabla h)(F_t,F_u,F_v)
  &=(\tfrac{\gamma_1'}{\gamma_1}-\tfrac{\alpha'}{\alpha}-2
  \tfrac{\gamma_2'}{\gamma_2}) h(F_u,F_v),\\
  (\nabla h)(F_t,F_v,F_v)
  &=(\tfrac{\gamma_1'}{\gamma_1}-\tfrac{\alpha'}{\alpha}-2
  \tfrac{\gamma_2'}{\gamma_2}) h(F_v,F_v),\\
  (\nabla h)(F_u,F_t,F_t)&=0=(\nabla h)(F_v,F_t,F_t),
\endalign$$
implying that $K_{F_t}$ restricted to the space spanned by $F_u$
and $F_v$ is a multiple of the identity.  Using the symmetries of
$K$ it now follows immediately that $F$ admits an $Z_3$-symmetry
or an $SO(2)$- symmetry.\qed
 \enddemo

 \specialhead 4.2.2 The second case: $\nu \equiv 0$ and $\l\neq \eta$ on $M^3$\endspecialhead

 Next, we consider the case that $b =\l^2-\eta^2$ and $\eta \neq
\l$ on $M^3$. Since by (4.35) and (4.34) $e_1(\eta -\l)=-\eta^2 -3
\l^2 -b+4\l \eta=4 \l (\eta -\l)$ we see that $\eta \neq \l$
everywhere on $M^3$ or nowhere.

We already have seen that $M^3$ admits a warped product structure.
The map $\varPhi$ we have constructed in Lemma 4.9 will not define an
immersion (cf. (4.40) and (4.41)). Anyhow, for a fixed point
$t_0$, we get from (4.26) - (4.29), (4.40) and (4.41), using the
notation $\tilde{\xi}=-(\eta+\l) e_1 + \xi$:
 $$\align
&D_{e_2}e_2 =\f_{22}^3 e_3 +\mu e_2 +\tilde{\xi},\\
&D_{e_2}e_3 =-\f_{22}^3 e_2 -\mu e_3,\\
&D_{e_3}e_2 =\f_{33}^2 e_3 -\mu e_3,\\
&D_{e_3}e_3 =\f_{33}^2 e_2 -\mu e_2 +\tilde{\xi},\\
&D_{e_i}\tilde{\xi}=0, \qquad i=2,3.
\endalign$$
 Thus, if $u$ and $v$ are local coordinates which span the second
distribution, then we can interpret $F(t_0,u,v)$ as an improper
affine sphere in a $3$-dimensional linear subspace.

Moreover, we see that this improper affine sphere is a paraboloid
provided that $\mu$ at $t_0$ vanishes identically (as a function
of $u$ and $v$). From the differential equations (4.20)
determining $\mu$, we see that this is the case exactly when $\mu$
vanishes identically, i.e.  when $M$ admits a pointwise
$SO(2)$-symmetry.

After applying a translation and a change of coordinates, we may
assume that
$$F(t_0,u,v)=(u,v,f(u,v),0),$$
with affine normal $\tilde{\xi}(t_0,u,v)=(0,0,1,0)$. To obtain
$e_1$ at $t_0$, we consider (4.24) and (4.25) and get that
 $$D_{e_i}(e_1 -(\eta -\l) F) = 0,\qquad i=2,3.$$
 Evaluating at $t=t_0$, this means that there exists a constant
vector $C$ such that $e_1(t_0,u,v)=(\eta -\l)(t_0) F(t_0,u,v) +C$.
Since $\eta\neq \l$ everywhere, we can write:
$$e_1(t_0,u,v)=\alpha_1 (u,v, f(u,v),\alpha_2),\tag4.43$$
where $\alpha_1\neq 0$ and we applied an equiaffine transformation
so that $C=(0,0,0,\alpha_1\alpha_2)$. To obtain information about
$\pt e_1$ we have that $D_{e_1}e_1= 2\l e_1 +\xi$ (cf. (4.21)) and
$\xi=\tilde{\xi} + (\eta+\l)e_1$ by the definition of
$\tilde{\xi}$. Also we know that $\tilde{\xi}(t_0,u,v)=(0,0,1,0)$
and by (4.39) - (4.41) that $D_{e_i}(\tau\tilde{\xi})=0$,
$i=1,2,3$. Taking suitable initial conditions for the function
$\tau$ ($\tau(t_0)=1$), we get that $\tau\tilde{\xi}=(0,0,1,0)$
and finally the following vector valued differential equation:
$$\pt e_1 =(\eta +3\l) e_1 +\tau^{-1} (0,0,1,0).\tag 4.44$$
Solving this differential equation, taking into account the
initial conditions (4.43) at $t=t_0$, we get that there exist
functions $\delta_1$ and $\delta_2$ depending only on $t$ such
that
$$e_1(t,u,v)= (\delta_1(t) u,\delta_1(t) v, \delta_1(t) (f(u,v)
  +\delta_2(t)), \alpha_2 \delta_1(t)),$$
  where $\delta_1(t_0)=\alpha_1$, $\delta_2(t_0)=0$, $\delta_1'(t)
=(\eta +3 \l) \delta_1(t)$ and $\delta_2'(t) =\delta_1^{-1}(t)
\tau^{-1}(t)$.  As $e_1(t,u,v) =\tfrac{\partial F}{\partial
  t}(t,u,v)$ and $F(t_0,u,v) =(u,v,f(u,v),0)$ it follows by
integration that
$$F(t,u,v)= (\gamma_1(t) u, \gamma_1(t) v, \gamma_1(t) f(u,v)
+\gamma_2(t) , \alpha_2 (\gamma_1(t)-1)),$$ where $\gamma_1'(t)
=\delta_1(t)$, $\gamma_1(t_0)=1$, $\gamma_2(t_0)=0$ and
$\gamma_2'(t) =\delta_1(t)\delta_2(t)$.  After applying an affine
transformation we have shown:
 \proclaim{Theorem 4.3}
 Let $M^3$ be a positive definite affine hypersurface of 
  $\R^4$ which admits a pointwise $SO(2)$- or $Z_3$-symmetry and with
  the globally defined functions satisfying $b
  =-\eta^2 +\l^2$ but $b \not\equiv 0$ on $M^3$. Then $M^3$ is affine
  equivalent with
  $$F:I\times N^2\to \R^4:(t,u,v)\mapsto (\gamma_1(t) u,
  \gamma_1(t) v, \gamma_1(t) f(u,v) +\gamma_2(t),\gamma_1(t)),$$
  where
  $\psi: N^2 \to \R^3:(u,v) \mapsto (u,v,f(u,v))$ is an
  improper affine sphere with affine normal $(0,0,1)$ and $\gamma:I\to
  \R^2$ is a curve.\hfill \newline Moreover, if $M^3$
  admits a pointwise $SO(2)$-symmetry then $N^2$ is an elliptic
  paraboloid.
 \endproclaim

 In the next theorem we again deal with the converse.

  \proclaim{Theorem 4.4}
 Let $\psi:N^2 \to \R^3:(u,v)\mapsto (u,v,f(u,v))$ be an
  improper affine sphere with affine normal
  $(0,0,1)$. Let $\gamma: I \to \R^2$ be a curve such that 
$$F(t,u,v)=(\gamma_1(t)u, \gamma_1(t)v,\gamma_1(t)f(u,v)
  +\gamma_2(t) ,\gamma_1(t)),$$
  defines a positive definite affine hypersurface. Then $F$
  admits a pointwise $Z_3$- or $SO(2)$-symmetry.
 \endproclaim
 \demo{Proof}
 We have
 $$\align
  &F_t=(\gamma_1' u,\gamma_1' v ,\gamma_1' f(u,v) +\gamma_2',\gamma_1'),\\
  &F_u=(\gamma_1,0,\gamma_1 f_u,0),\\
  &F_v=(0,\gamma_1,\gamma_1 f_v,0),\\
  &F_{tt}=(\gamma_1'' u, \gamma_1''v ,\gamma_1'' f(u,v) +
  \gamma_2'',\gamma_1'')=\tfrac{\gamma_1''}{\gamma_1'}F_t
  +\tfrac{(\gamma_2''\gamma_1'-
    \gamma_1''\gamma_2')}{\gamma_1'}(0,0,1,0),\\
  &F_{ut}=\tfrac{\gamma_1'}{\gamma_1} F_u,\\
  &F_{vt}=\tfrac{\gamma_1'}{\gamma_1} F_v,\\
  &F_{uu}=(0,0,f_{uu} \gamma_1,0),\\
  &F_{uv}=(0,0,f_{uv} \gamma_1,0),\\
  &F_{vv}=(0,0,f_{vv}\gamma_1,0).
\endalign$$
This implies that $F$ defines a nondegenerate affine immersion
provided that $\gamma_1 \gamma_1'(\gamma_2''\gamma_1'-
\gamma_1''\gamma_2') \neq 0$.  We moreover see that this immersion
is definite provided that the improper affine sphere is positive
definite and $\gamma_1\gamma_1'(\gamma_2''\gamma_1'-
\gamma_1''\gamma_2')>0$ or when the improper affine sphere is
negative definite and $\gamma_1\gamma_1'(\gamma_2''\gamma_1'-
\gamma_1''\gamma_2')<0$.  As the proof in both cases is similar,
we will only treat the first case here. It easily follows that we
can write the affine normal $\xi$ as:
$$
\xi=\alpha(t) (0,0,1,0) +\beta(t) F_t,
$$
where $(\gamma_2''\gamma_1'-\gamma_1'' \gamma_2')=
\gamma_1^2(\gamma_1')^3 \alpha^5$ and $\alpha'
+\tfrac{\beta(\gamma_2''\gamma_1'-\gamma_1''\gamma_2')}{\gamma_1'}=0$.
Taking $e_1$ in the direction of $F_t$, we see that $F_u$ and
$F_v$ are orthogonal to $e_1$. It is also clear that $S$
restricted to the space spanned by $F_u$ and $F_v$ is a multiple
of the identity, and $S(F_t)=a F_t$, since $S$ is symmetric.
Moreover, we have that
 $$\align
  (\nabla h)(F_t,F_u,F_u)
  &=(-\tfrac{\gamma_1'}{\gamma_1}-\tfrac{\alpha'}{\alpha})
  h(F_u,F_u),\\
  (\nabla h)(F_t,F_u,F_v)
  &=(-\tfrac{\gamma_1'}{\gamma_1}-\tfrac{\alpha'}{\alpha})
  h(F_u,F_v),\\
  (\nabla h)(F_t,F_v,F_v)
  &=(-\tfrac{\gamma_1'}{\gamma_1}-\tfrac{\alpha'}{\alpha})
  h(F_v,F_v),\\
  (\nabla h)(F_u,F_t,F_t)&=0=(\nabla h)(F_v,F_t,F_t),
\endalign$$
implying that $K_{F_t}$ restricted to the space spanned by $F_u$
and $F_v$ is a multiple of the identity.  Using the symmetries of
$K$ it now follows immediately that $F$ admits an $Z_3$-symmetry
or an $SO(2)$- symmetry.\qed
 \enddemo

 \specialhead 4.2.3 The third case: $\nu \equiv 0$ and $\l= \eta$ on
 $M^3$\endspecialhead

 The final case now is that $b =\l^2-\eta^2$ and $\eta = \l$ on the
whole of $M^3$, i.~e. $b=0$. This is dealt with in the following theorem:
\proclaim{Theorem 4.5}
 Let $M^3$ be a positive definite hypersurface of $\R^4$ which admits
  a pointwise $SO(2)$- or $Z_3$-symmetry and with the globally defined
  functions satisfying $b =0$ and $\eta =\l$ on $M^3$. Then $M^3$ is
  affine equivalent to
  $$F:I\times N^2\to \R^4:(t,u,v)\mapsto (u, v, f(u,v)
  +\gamma_2(t),\gamma_1(t)),$$
  where $\psi: N^2 \to \R^3:(u,v)
  \mapsto (u,v,f(u,v))$ is an improper affine sphere with affine
  normal $(0,0,1)$ and $\gamma:I\to \R^2$ is a curve.\hfill 
\newline Moreover, if $M^3$ admits a pointwise
  $SO(2)$-symmetry then $N^2$ is an elliptic paraboloid.
 \endproclaim
  \demo{Proof}
  We proceed in the same way as in Theorem 4.3.
  We again use that $M^3$ admits a warped product structure and we fix a
  parameter $t_0$. At the point $t_0$, we have for $\tilde{\xi}=-(\eta+\l)
e_1 + \xi = -2\l e_1 + \xi$:
 $$
\align
&D_{e_2}e_2 =\f_{22}^3 e_3 +\mu e_2 +\tilde{\xi},\\
&D_{e_2}e_3 =-\f_{22}^3 e_2 -\mu e_3,\\
&D_{e_3}e_2 =\f_{33}^2 e_3 -\mu e_3,\\
&D_{e_3}e_3 =\f_{33}^2 e_2 -\mu e_2 +\tilde{\xi},\\
&D_{e_i}\tilde{\xi}=0, \qquad i=2,3.
\endalign$$

Thus, if $u$ and $v$ are local coordinates which span the second
distribution, then we can interpret $F(t_0,u,v)$ as an improper
affine sphere in a $3$-dimensional linear subspace.

Moreover, we see that this improper affine sphere is a paraboloid
provided that $\mu$ at $t_0$ vanishes identically (as a function
of $u$ and $v$). From the differential equations (4.20)
determining $\mu$, we see that this is the case exactly when $\mu$
vanishes identically, i.e.  when $M$ admits a pointwise
$SO(2)$-symmetry.

After applying a translation and a change of coordinates, we may
assume that
$$
F(t_0,u,v)=(u,v,f(u,v),0),
$$
with affine normal $\tilde{\xi}(t_0,u,v)=(0,0,1,0)$. To obtain
$e_1$ at $t_0$, we consider (4.24) and (4.25) and get that
$$
D_{e_i}e_1 =(\eta -\l) e_i = 0,\qquad i=2,3.
$$
It follows that $e_1(t_0,u,v)$ is a constant vector field. As it
is transversal, we may assume that there exists an $\alpha\neq 0$
such that
$$
e_1(t_0,u,v)= (0,0,0,\alpha).
$$
As $e_1$ is determined by the differential equation (cf. (4.44)):
$$
\frac{\partial e_1}{\partial t} =4\l  e_1 +\tau^{-1} (0,0,1,0),
$$
it follows that
$$
e_1(t,u,v)= (0,0,\delta_2(t),\delta_1(t)),
$$
where $\delta_2(t_0)=0$, $\delta_2'(t)=4\l(t) \delta_2(t)+
\tau^{-1}(t)$, $\delta_1(t_0)=\alpha$ and $\delta_1'(t)=4\l(t)
\delta_1(t)$. Integrating once more with respect to $t$ we obtain
that
$$
F(t,u,v)= (u,v, f(u,v)+\gamma_2(t),\gamma_1(t)),
$$
for some functions $\gamma_1$ and $\gamma_2$ with
$\gamma_i'=\delta_i$ and $\gamma_i(t_0)=0$ for $i=1,2$.\qed
  \enddemo

In the next theorem we deal with the converse.

\proclaim{Theorem 4.6}
 Let $\psi:N^2 \to \R^3:(u,v)\mapsto (u,v,f(u,v))$ be an
  improper affine sphere with affine normal $(0,0,1)$. Let
  $\gamma: I \to \R^2$ be a curve such that
  $$F(t,u,v)=(u, v,f(u,v) +\gamma_2(t) ,\gamma_1(t))$$
  defines a positive definite affine hypersurface. Then $F$ 
  admits a pointwise $Z_3$- or $SO(2)$-symmetry.
 \endproclaim
  \demo{Proof}
  We have
  $$\align
  &F_t=(0,0,\gamma_2',\gamma_1'),\\
  &F_u=(1,0,f_u,0),\\
  &F_v=(0,1,f_v,0),\\
  &F_{tt}=(0,0,\gamma_2'',\gamma_1'')=\tfrac{\gamma_1''}{\gamma_1'}F_t
  +\tfrac{(\gamma_2''\gamma_1'-
    \gamma_1''\gamma_2')}{\gamma_1'}(0,0,1,0),\\
  &F_{ut}=F_{vt}=0,\\
  &F_{uu}=(0,0,f_{uu},0),\\
  &F_{uv}=(0,0,f_{uv},0),\\
  &F_{vv}=(0,0,f_{vv},0).
\endalign$$
This implies that $F$ defines a nondegenerate affine immersion
provided that $\gamma_1'(\gamma_2''\gamma_1'-\gamma_1''\gamma_2')
\neq 0$.  We moreover see that this immersion is definite provided
that the improper affine sphere is positive definite and
$(\gamma_2''\gamma_1'-\gamma_1''\gamma_2')\gamma_1'>0$ or when the
improper affine sphere is negative definite and
$(\gamma_2''\gamma_1'-\gamma_1''\gamma_2')\gamma_1'<0$.  As the
proof in both cases is similar, we will only treat the first case
here. It easily follows that we can write the affine normal $\xi$
as:
$$
\xi=\alpha(t) (0,0,1,0) +\beta(t) F_t,
$$
where $(\gamma_2''\gamma_1'-\gamma_1''\gamma_2') = (\gamma_1')^3
\alpha^5$ and $\alpha'
+\tfrac{\beta(\gamma_2''\gamma_1'-\gamma_1''\gamma_2')}{\gamma_1'}=0$.
Taking now $e_1$ in the direction of $F_t$, we see that $F_u$ and
$F_v$ are orthogonal to $e_1$. It is also clear that $SF_u =
SF_v=0$, and $S(F_t)=a F_t$, since $S$ is symmetric. Moreover, we
have that $$\align
  (\nabla h)(F_t,F_u,F_u) &=-\tfrac{\alpha'}{\alpha} h(F_u,F_u)\\
  (\nabla h)(F_t,F_u,F_v) &=-\tfrac{\alpha'}{\alpha} h(F_u,F_v)\\
  (\nabla h)(F_t,F_v,F_v) &=-\tfrac{\alpha'}{\alpha} h(F_v,F_v),\\
  (\nabla h)(F_u,F_t,F_t)&=0=(\nabla h)(F_v,F_t,F_t),
\endalign$$
implying that $K_{F_t}$ restricted to the space spanned by $F_u$
and $F_v$ is a multiple of the identity.  Using the symmetries of
$K$ it now follows immediately that $F$ admits an $Z_3$-symmetry
or an $SO(2)$-symmetry.  \qed
  \enddemo

 \enddocument
 \end